# THE CONTINUUM OF THE SURREAL NUMBERS REVISITED. THE SURREAL NUMBERS DEFINED THROUGH TRANSFINITE FUNDAMENTAL CAUCHY SEQUENCES

## Constantine E. Kyritsis*

*\* Associate Prof. of the University of Ioannina, Greece.*
*ckiritsi@uoi.gr  C_kyrisis@yahoo.com, Dept. Accounting-Finance*
*Psathaki Preveza 48100*

## Abstract

In this treatise on the theory of the continuum of the surreal numbers of J.H. Conway, is proved ,that the three different techniques and hierarchies of the continuums of the transfinite real numbers of Glayzal A. (1937) defined through transfinite power series , of the surreal numbers of J.H. Conway (1976) defined by Dedekind cuts ,and of the ordinal real numbers of K. E. Kyritsis (1992) defined by fundamental Cauchy transfinite sequences, give by inductive limit or union the same class and continuum of infinite numbers. This is quite remarkable and is the analogue in the transfinite numbers, of that the real numbers can be constructed either as decimal power series, or by Dedekind cuts, or by Cauchy fundamental sequences.

**Subject Classification of AMS** 03,04,08,13,46

## PROLOGUE

In 1976 J.H. Conway published a book with title "On Numbers and games", where utilizing the very well-known technique o the Dedekind cuts and applied not only on the finite natural and rational numbers on the ordinal numbers as well, a new continuum is created analogous to the real numbers but vastly larger and finer, that he called **the surreal numbers.** In this treatise we re-create the surreal number with the other very well known technique of the Cauchy fundamental sequences. Actually a new vastly large and fine continuum is also created that the author called **ordinal real numbers**. This creative work was carried out during 1990-1991 in the birth Greek island of the Pythagoras called Samos, and in the mathematical Department of the University of the Aegean. When this was complete, the author discovered that still a third older author A, Glayzal in 1937, had created also such a vast and fine continuum through transfinite power series that he called **transfinite real numbers**. Then in 1992 the author proved that the three continuums of transfinite real numbers, of the surreal numbers and of the ordinal real numbers is one and the same continuum. This fact, includes the very well-known fact in the foundations of the mathematics that the usual real numbers can be created either by decimal power series, or by Dedekind cuts or by Cauchy fundamental sequences. This work was partitioned in 5 separate papers. These papers remained unpublished as the author had an Odyssey of different jobs, and eventually were published during 2017, in the proceedings https://books.google.gr/books?id=BSUsDwAAQBAJ&pg **pp 233-292**

 of the conference 1st INTERNATIONAL CONFERENCE ON QUANTITATIVE, SOCIAL, BIOMEDICAL AND ECONOMIC ISSUES 29-30 JUNE 2017 https://www.linkedin.com/pulse/1st-international-conference-quantitative-social-economic-frangos/ .

http://icqsbei2017.weebly.com/

The next   table summarizes the equivalence of the hree techniques so as to create this continuum of the surreal numbers.

Table 1

| UNIFICATION | HOW IT IS CRERATED | THE RESULT |
|---|---|---|
| **TRANFINITE REAL NUMBERS** by A. Glayzal  1937 | **Transfinite power series** | The  continuum No of the surreal numbers |
| **THE SURREAL NUMBERS** by J.H. Conway 1976 | **Dedekind cuts** | The  continuum No of the surreal numbers |
| **THE ORDINAL REAL NUMBERS** by K.E. Kyritsis 1991 | **Cauchy fundamental sequences** | The  continuum No of the surreal numbers |

We should not understand that with the current theory we suggest direct applications in the physical sciences.  Not at all!  **Matter is always finite**. Actually not even the real numbers are fully appropriate for the physical reality because they are based on the infinite too which does not exist in the material reality. This has been described in more detailed by the famous Nobel prize winner physicist E. Schrödinger in his book "Science and Humanity" (see [ Schrödinger E. 1961]. That is why the author has developed the digital or natural real numbers without the infinite with the corresponding Euclidean geometry and also Differential and Integral calculus, which is logically different from the classical. (See [Kyritsis 2017] and [Kyritsis 2019]) But the ordinal numbers and the surreal numbers reflect more the human **consciousness and perceptions** rather than properties of the physical material reality. Still such a discipline as the study of the continuum of the surreal numbers is an excellent spiritual, mental and metaphysical meditative practice probably better than many other metaphysical spiritual systems. ***It is certainly an active reminding to the scientists that the ontology of the universe is not only the finite matter but also the infinite perceptive consciousness.***

**Alternative algebraic definitions of the Hessenberg   natural operations   in the ordinal numbers. The ORDINAL NATURAL NUMBERS 1**

**Constantine E. Kyritsis***

*\* Associate Prof. of the University of Ioannina, Greece.* ckiritsi@teiep.gr   C_kyrisis@yahoo.com, *Dept. Accounting-Finance Psathaki Preveza 48100*

## Abstract

This paper proves prerequisite results for the theory of Ordinal Real Numbers. In this paper, is proved that any field-inherited abelian operations and the Hessenberg operations ,in the ordinal numbers coincide. It is given an algebraic characterization of the Hessenberg operations ,that can be described as an abelian, well- ordered, double monoid with cancelation laws.

**§0.  Introduction**  This is a series of five papers that have as goal the definition of topological complete linearly ordered fields (continuum of transfinite numbers) that include the real numbers and are obtained from the ordinal numbers in a method analogous to the way that Cauchy derived the real numbers from the natural numbers. We shall also prove that this continuum is nothing lase than the **continuum of the surreal numbers** of J.H. Conway.

We should not understand that with the current theory we suggest direct applications in the physical sciences. Not at all! **Matter is always finite**. Actually not even the real numbers are fully appropriate for the physical reality because they are based on the infinite too which does not exist in the material reality. This has been described in more detailed by the famous Nobel prize winner physicist E. Schrödinger in his book "Science and Humanity" (see [ Schrödinger E.  1961]. That is why the author has developed the digital or natural real numbers without the infinite with the corresponding Euclidean geometry and also Differential and Integral calculus, which is logically different from the classical. (See [Kyritsis 2017] and [Kyritsis 2019]) But the ordinal numbers and the surreal numbers reflect more the human **consciousness and perceptions** rather than properties of the physical material reality. Still such a discipline as the study of the continuum of the surreal numbers is an excellent spiritual, mental and metaphysical meditative practice probably better than many other metaphysical spiritual systems. *It is certainly an active reminding to the scientists that the ontology of the universe is not only the finite matter but also the infinite perceptive consciousness.*

The author started and completed this research in the island of Samos during 1990-1992 .

In this paper are studied the Hessenberg operations in the ordinal numbers, from an algebraic point of view.

The main results are the characterization theorems 9,10 . They are characterizations of the Hessenberg operations as

a) field-inherited operations in the ordinal numbers, that satisfy two inductive properties .(see proposition 10)

b) operations that satisfy a number of purely algebraic properties, that could be called in short operations of a well-ordered commutative semiring with unit ;(see Lemma 0, proposition 9).

In a next paper I shall give two more algebraic characterisations of the Hessenberg operations as

c) operations defined by transfinite induction in the ordinal numbers and by two recursive rules ,

d) operations of the free semirings in the category of abelian semirings; or as the operations of the formal polynomial algebras of the category of abelian semirings. These characterisations of the Hessenberg operations are independed of the standard non-commutative operations of the ordinal numbers and can be considered as alternative and simpler definitions of them (especially the c),d)).

In particular it is proved that the Hessenberg natural operations are free finitary operations ; We make use of rudimentary techniques relevant to K-theory and Universal Algebra .

The main application of the present results is in the definition of the **ordinal real numbers**.(see [ Kyritsis C.E.1991] .By making use of the present results and techniques it is proved in [Kyritsis C.E.1991] that all the three techniques and Hierarchies of *transfinite real numbers* ,see [ Glayzal A 1937],of *surreal numbers*, see [ Conway J.H. 1976], of *ordinal*

*real numbers* see [ Kyritsis C.E.1991], give by inductive limit or union, the same class of numbers ,already known as the class No. We refer to the class No as **"the totally ordered Newton-Leibniz realm of numbers"**).

## 1. Two algebraic characterizations of the Hessenberg operations in the ordinal numbers.

Let us denote by F a linearly ordered field of characteristic $\omega$ (also said of characteristic 0).Let us denote by h a mapping of an initial segment of ordinal numbers, denoted by W(a),in F such that it is 1-1, order preserving and h(0)=0, h(1)=1, h(s(b))=b+1,where s(b) is the sequent of the b , b<a, the b+1 is in the field operations the set h(W(a)) is closed in the field addition and multiplication . We shall call field-inherited operations in the ordinal numbers of W(a), the operations induced by the field, in the initial segment W(a).

(For a reference to standard symbolisms and definitions for ordinal numbers, see [Cohn P.M. 1965] p.1-36 also [ Kutatowski K.-Mostowski A. 1968]).

The following properties hold for these field-inherited abelian operations for the ordinals of w(a) ( in that case, it is needless to say that a is a limit ordinal).

**Lemma 1.** For the field-inherited abelian operations in the initial segment w(a) of ordinal numbers, the followings hold (x,y,z,c,x',y', ∈w(a)).

0)  s(x) = x+1              for every x ∈ w(a)

1)  x+y = y+x                ,   x.y = y.x

2)  x+(y+z) = (x+y)+z          ,   x.(y.z) = (x.y).z

3)  x+0 = 0+x = x      ,

    x.0 = 0

    x.1 = 1.x = x

4)  x(y+z) = xy+xz

5)  $\left.\begin{array}{ll} x+c = y+c \Rightarrow x = y & x \cdot c = y \cdot c \\ c \in w(a) & c \neq 0 \end{array}\right\} \Rightarrow x = y$

6) If x>y, x'>y' then x+x'>y+y', and xx'+yy'>xy'+yx'

The proof of the previous lemma is direct from the properties of a linearly ordered field.

We mention two more properties that they will be of significance in the followings paragraphs.

7) The w(x+y) is a cofinal set with the $\{w(x)+x \cup x+w(y)\}$ and we write $cf(w(x+y)) = cf\{w(x)+y \cup x+w(y)\}$.

8) The W(x.y) is cofinal set with the

$h^{-1}(\{h(y)h(w(x))+h(x)h(w(y))-h(w(x)).h(w(y))\})$

and we write

$cf(w(x.y)) = cf \; h^{-1} \; (\{h(y)h(w(x))+h(x)h(w(y))-h(w(x)).h(w(y))\})$

To continue our argument we need a many-variables form of transfinite induction.

Let $a_i$ $i = 1,...,n$ $n \in N$ ordinal numbers and $(b_1,...b_n) \in w(\alpha_i)x...xw(\alpha_n)$. We define as simultenous initial segment of n-variables defined by $(b_1,...b_n)$, the set $w((b_1,...b_n)) = \cup \; w(b_1)x...x\{b_1\}x...xw(b_n)$ for every

$i \in I \subseteq \{1,...,n\}$

$I \underset{\neq}{\subseteq} \{1,...n\}$, or $w((b_1,...,b_n))= \cup w\{b_1\}x...xw\{b_i\}x...x\{b_n\}$ for

$i \in I \subseteq \{1,...,n\}$

every $I \subseteq \{1,...n\}$ with $I \neq \varnothing$.

**Lemma 2. (many-variables transfinite induction)**

Let $A \subseteq w(a_1)x...xw(a_n)$ such that

1. $(0,...,0) \in A$

2. For every $(b_1,...b_n) \in w(a_1)x...xw(a_n)$ it holds that $w((b_1,...,b_n)) \subseteq A$ $\Rightarrow (b_1,...,b_n) \in A$.

By 1.2. we infer that $A = w(a_1)x...xw(a_n)$.

**Lemma 3. (many-variables definition by transfinite induction)**

Let a set $A$ and ordinal numberss $a_1,...,a_n$. Let a set denoted by $B$, such that it is sufficient for inductive rules $h\ B \rightarrow A$ :

This means that :

a) $B \underset{(b_1,...,b_n)}{\subseteq} \cup A^{w((b_1,...,b_n))},\ (b_1,...,b_n) \in w(a_1)x...xw(\alpha_n)$

The set $B$ is a set of functions, denoted by $f_{(b_1,...,b_n)}$ and defined on simulteneous initial segments with values in $A$. $f_{(b_1,...,b_n)} : w((b_1,...,b_n)) \rightarrow A$.

If $f_{(b_1,...,b_n)} \in B$ and $c_1 < b_1,...,c_n < b_n$ then $\dfrac{f(b_1,...,b_n)}{w((c_1,...c_n))} \in B$

b) For every $(c_1,...,c_n) \in w(\alpha_1)x...xw(\alpha_n)$ there is a $f_{(c_1,...,c_n)}$ such that $f_{(c_1,...,c_n)} \in B$

c) Let $f \in B^{w((\alpha_1,...,\alpha_n))}$ and let us denote the value of $f$ at $(b_1,...,b_n)$ with $f_{(b_1,...,b_n)}$. Let us suppose that it holds that whenever $c_1 \le b_1, b'_1,...,c_n \le b_n, b'_n$, $(b_1,...,b_n),(b'_1,...,b'_n) \in$

$w(\alpha_1)x...xw(\alpha_n)$ then $f_{(b_1,...,b_n)} \Big/ w_{((c_1,...c_n))} = f_{(b'_1,...,b'_n)_{w((c_1,...c_n))}}$

Then let us suppose that we get as a consequence that the function defined by $g(b_1,...,b_n) = f_{(b_1,...,b_n)}(b_1,...,b_n),(b_1,...,b_n) \in w((\alpha_1,...,\alpha_n))$, belongs to $B$

It holds that for every function $h:B \rightarrow A$

(called many -variables transfinite inductive rule ) there is one and only one function $f$ defined on $w(a_1)x...xw(a_n)$ with values in $A$ such that for every $(b_1,...,b_n) \in w(\alpha_1)x...xw(\alpha_n)$ it holds that $f(b_1,...,b_n) = h\Big(f\Big/w((b_1,...,b_n))\Big)$

**Remark**:We notice that even for one variable this version of the definition by transfinite induction is somehow different from that which appears usually in the bibliography (e.g.see [Kutayowski K. –Mostowski A 1968] §4 pp 233 ).It uses not all the set $A^{w(\alpha)}$ ,but only a subset of it,  sufficient for recursive rules .The proof ,for one variable, is nevertheless exactly the same as with the ordinary version .

In order to save space and because the proofs are  not directly relevant  to our subject  we will not give the  proofs of lemma 2  and 3 but we will mention that they are  analogous, without  serious difficulties, to the ones with  one-variable only  (see  e.g. [ Kutayowski K. –Mostowski A 1968],[ Lang S.1984]).

## Proposition 4. (Uniqueness)

Any two pairs of field-inherited operations in  the initial segment w(a) of ordinals, satisfying properties  7,8 of lemma 1 (a  is a limit ordinal ) are isomorphic .

**Proof**: Let a monomorphic embedding denoted by h of w(a), as is described in  the beggining of the paragraph in two linearly ordred fields denoted by   $F_1$, $F_2$.

Let the two pairs of field inherited operations in    w(a) be denoted by $((+,.)\,(^{\oplus,\circ})$ respectively.   They    satisfy the properties 0.1.2.3.4.5.6.7.8.of lemma 1 . Suppose that the  operations +, $\oplus$ coincide for   the  set $w((b_1,b_2))$ $\subseteq w(\alpha)^2$ where

$b_1,b_2 \in w(\alpha)$   $w((b_1,b_2)) = w(b_1)xw(b_2) \cup (\{b_1\}xw(b_2)) \cup (w(b_1)x\{b_2\})$ . Then  by property   7   $b_1 + b_2 = s(b_1 + w(b_2) \cup w(b_1) + b_2) =$ (by  the   hypothesis of transfinite induction) =

$$s(b_1 \oplus w(b_2) \cup w(b_1) \oplus b_2) = b_1 \oplus b_2;$$

Where by S(A) we denote the sequent of the set A .Thus by lemma 2 the operations +, $\oplus$ coincide on  w(a)xw(a).

Then the set w(a) is an ordered abelian  monoid  relative  to addition, with cancelation law.

The Grothendieck groups of w(a) for both +, and   $\oplus$ coincide, and we denote it by k(w(a)) (see for  the  definition  of Grothendieck group [Lang S. 1984] Ch1 §4 p. 44  or  [Cohn P.M. 1965]  ch  vii  §3 pp 263  ). Thus  also  the opposite -x of an element x  of w(a) is  the  same  in  the

Grothendieck group k(w(a)) of w(a) for both the two operations + and $\oplus$ .

Suppose also that the operations $\overset{\bullet}{\phantom{.}},\circ$ coincide for the set w((b$_1$, b$_2$)). Then by property 8

$$b_1 \bullet b_2 = s\left(h^{-1}\left(\left\{h(b_1)h(w(b_2)) + h(w(b_1)) \bullet h(b_2) - h(w(b_1)) \bullet h(w(b_2))\right\}\right)\right) =$$

(because +, and $\oplus$ are isomorphic and the hypothesis of transfinite induction for

$\bullet,\circ) = s\left(h^{-1}\left(\left\{h(b_1) \circ h(w(b_2)) \oplus h(w(b_1)) \circ h(b_2) - h(w(b_1)) \circ h(w(b_2))\right\}\right)\right) = b_1 \circ b_2$ . Hence by lemma 2 the two operations $\overset{\bullet}{\phantom{.}},\circ$ coincide on the whole set w(a)xw(a) Q.E.D.

The next step is to find the relation of field-inherited operations in an initial segment of ordinals with the Hessenberg operations . It will turn out that, if they satisfy the properties 7.8.of lemma 1 ,then they are nothing else than the Hessenberg-Conway natural sum and product (see [Kutatoski K.Mostowski A. 1968 ] ch VII §7 p. 252-253 exercises 1. 2. 3.) and [ Frankel A.A.1953] pp. 591-594 also [Conway J.H. 1976] ch2 p. 27-28).

The way in which the Hessenberg operations are defined, traditionally , depends on the standard non-commutative operation on ordinals.

In order to define the Hessenberg-Conway operations in the traditional way ,we remind that :

## Lemma 6 (Cantor normal form).

For every ordinal a there exists a natural number n and finite sequences : b$_1$,...,b$_2$ of natural numbers and ordinal numbers a$_1$,a$_2$,...a$_n$ with a$_1$>...>a$_n$ such that $\alpha = \omega^{a_1} b_1 + ... + \omega^{a_n} b_n$ (For a proof see for instance [Kutatowski K.-Mostowski A. 1968] ch VII §7 p. 248-251).

Then we get for the two ordinal numbers  $\alpha$, $b$, by adding terms with zero coefficients, to make  their  Cantor normal  forms of  equal   length ,
that      $\alpha = \omega^{\xi_1}p_1 + ... + \omega^{\xi_n}p_n \quad b = \omega^{\xi_1}q_1 + ... \omega^{\xi_n}q_n$;

we  define  the  natural  sum  (we  denote  it  by  (+))  with
$\alpha(+)b = \omega^{\xi_1}\left(p_1 + q_1\right) + ... + \omega^{\xi_n}\left(p_n + q_n\right)$. The  natural   product, denoted by $\alpha(.)b$
is defined to be  the  ordinal  arising  by multiplication (using distributive
and associative laws) from the Cantor normal forms of $a$ and $b$ and  by
using  the  rule: $\omega^x(.)\omega^y = \omega^{x(+)y}$ to multiply powers of $\omega$. As a  result  we
get for instance that

## Remark 7

1)  The normal form of $a$ can also be written in the standard Hessenberg-Conway operations that is

$$\alpha = \omega^{\xi_1}(.)p_1(+)...(+)\omega^{\xi_n}(.)p_n.$$

2)  The sum $a(+)b$ is an increasing function of $a$ and $b$.

3)  If  $\xi < \omega^{\omega^\alpha}$ and  $\eta < \omega^{\omega^\alpha}$   then  $\xi.\eta < \omega^{\omega^\alpha}$ for ordinals  $\zeta$, $\eta$, $\alpha$ and
conversely if an ordinal $j$  satisfies the condition: "if $\xi < j$ and $\eta < j$ then $\xi.\eta < j$"
then  there  exists  an  ordinal number $\alpha$ such that  $\xi = \omega^{\omega^\alpha}$ ;we call ordinal
numbers of the        type  $\omega^{\omega^\alpha}$ **principal ordinals**  of the  Hessenberg
operations.(see [ Kutatowski K.-Mostowski A. 1968] ch vii paragraph 7,p
253) This  has also as  a consequence  that  we define the Hessenberg-
Conway          natural operations only for initial segments of the type W
$\left(\omega^{\omega^\alpha}\right)$ for  some  ordinal  number  $\alpha$ (we  will  call them   principal initial
segments ).

4)   The Hessenberg-Conway natural operation restricted on the set of
Natural numbers coincide with  the  ordinary  sum and product of natural
numbers.

5)  The  operation   "powers  of   $\omega$" ,  through   the Hessendberg-
Conway natural operation, can be  defined  as follows:

a)   $\omega^{(0)} = 1$ $\omega^{(1)} = \omega$ if $\xi$ is a limit ordinal $\omega^{(\xi)} = \sup \omega^{(\eta)}$

$$\eta < \xi$$

b)   If $\xi$ is not a limit ordinal then there exists an ordinal $\eta$ with $\eta(+)1 = s(\eta) = \xi$ and we define $\omega^{(\xi)} = \omega^{\eta}(.)\omega$.

It holds, (this happens especially for the base $\omega$), that these "natural powers" of $\omega$ coincide with the standard powers of $\omega$ defined through the standard non-commutative multiplication of ordinal numbers (this can be proved with transfinite induction since $\omega^{\eta}.\omega = \omega^{(\eta)} \circ \omega$. This gives us the right to express any ordinal number $\alpha$ in Cantor normal form, exclusively with natural operations:

$$\alpha = \omega^{(\xi_1)}(.)p_1(+)...(+)\omega^{(\xi_n)}(.)p_n$$

6)   Also we notice that, the natural difference denoted by a (-) b, of two ordinals a, b in Cantor normal forms $a = \omega^{(\xi_1)}(.)p_1(+)...(+)\omega^{(\xi_n)}(.)p_n$

$$b = \omega^{(\xi_1)}(.)q_1(+)...(+)\omega^{(\xi_n)}(.)q_n$$ , is defined only if $p_1 \geq q_1 ... p_n \geq q_n$.

7)   We notice that if $\xi_i < \xi_j$ for two ordinals then $\omega^{\xi_i} < \omega^{\xi_j}$ but also $\omega^{\xi_i}.a < \omega^{\xi_j}.b$ for every pair of non-zero natural numbers a,b. (in the standard non-commutative operations on ordinals). But this has as consequence that the ordering of a finite set of ordinal numbers in Cantor normal form (normalizing the Cantor normal forms by adding terms with zero coefficients so that all of them have the same set of exponents) is isomorphic (similar) to the lexicographical ordering of the coefficients of the normal forms.

**Proposition 8.** For every principal initial segment of ordinal numbers, the Hessenberg natural operations satisfy the properties 0.1.2.3.4.5.6.7.8. of lemma 1.

**Remark**. From the moment we have proved the properties 0.1.2.3.4.5.6. for the natural operations in the principal initial segment w(a), there is the Grothendieck group k(w(a)) of the monoid relative to sum, w(a) such that the w(a) is monomorphicaly embedded in k(w(a) (because of cancelation law) and also there is an ordering in k(w(a)) that restricted on w(a) coincides with the standard ordering in w(a).

Then the difference that occurs in property 8 has meaning and also the statement of property 8 itself has meaning (see [Lang S. 1984] Ch I §9

p. 44).We denote by h  the monomorphism of  the W(a)  in  the  K(W(a)) .

**Proof.** The  properties  0.1.2.3.4.  are  directly  proved  from  the definition of the natural operations.  Let us check  the property 5. Namely , the  cancelation  laws. Let us suppose  that  y,x,c  are ordinal numbers with y,x,c $\in w\left(\omega^{\omega^a}\right)$  for some  ordinal  a and their  Cantor  normal  form , in  natural  operations ,  are

$$x = \omega^{(\xi_1)}(.)p_1(+)...(+)\omega^{(\xi_n)}(.)p_n$$

$$c = \omega^{(\xi_1)}(.)c_1(+)...(+)\omega^{(\xi_n)}(.)c_n$$

$$y = \omega^{(\xi_1)}(.)q_1(+)...(+)\omega^{(\xi_n)}(.)q_n \quad p_i, c_i, y_1 \in No$$

then

$$x(+)c = \omega^{(\xi_1)}(.)\left(p_1 + c_1\right)(+)...(+)\omega^{(\xi_n)}(.)\left(p_n + c_n\right)$$

$$y(+)c = \omega^{(\xi_1)}(.)\left(q_1 + c_1\right)(+)...(+)\omega^{(\xi_n)}(.)\left(q_n + c_n\right)$$

hence

$$x(+)c=y(+)c \Rightarrow p_i +c_i = q_i +c_i \qquad i = 1,...,n$$

and by cancelation law for addition  in  natural  numbers  we deduce that $p_i=q_i$  i=1,...,n  hence x = y.

Also

$$x(.)c = \sum_{1\le i\ j\le n} \omega^{\xi_i(+)\xi_j}(.)\left(p_i.c_j\right)$$
and

$$y(.)c = \sum_{1\le i\ j\le n} \omega^{\xi_i(+)\xi_j}(.)\left(q_i.c_j\right)$$
and if $c \ne 0$

and $x(.)c = y(.)c$ then $p_i.c_j = q_i.c_j$ with  not  all  of  $c_j$  equal  to zero.  Say $c_{j_0} \ne 0$, then $p_i.c_{j_0} = q_i.c_{j_0}$ for  every  i = 1,...,n hence $p_i = q_i$ and x=y.

Let us check the  property  6. The first part of  property  6 is immediate from Remark 7, 2. Let, furthermore, x',y' $\subset w\left(\omega^{\omega^a}\right)$  with Cantor normal

form (changing the $\xi_i$, in order to have the same exponents for all x, y, x', y')

$$x' = \sum_1^n \omega^{(\xi_i)}(.)p_i' \qquad\qquad y' = \sum_1^n \omega^{(\xi_i)}(.)q_i'$$

with $p_i', q_i' \in$ No and with summation interpreted as natural sum. By hypothesis x'>y', x>y.

Then

$$x(.)x' = \sum_{ij} \omega^{\xi_i(+)\xi_j}(.)\left(p_i . p_j'\right) \qquad \text{and}$$

$$y(.)y' = \sum_{ij} \omega^{\xi_i(+)\xi_j}(.)\left(q_i . q_j'\right)$$

and the coefficient of the monomial of greatest exporent of x(.)x'(+)y(.)y' is $p_1 p_1' + q_1 q_1'$ and of x(.)y'(+)y(.)x' is $p_1 q_1' + q_1 p_1'$. But $p_1 p_1' + q_1 q_1' - p_1 q_1' - q_1 p_1' = p_1(p_1'-q_1') - q_1(p_1'-q_1') = (p_1-q_1).\ (p_1'-q_1') > 0$ which is a product of the positive factors $p_1-q_1$, $(p_1'-q_1')$ hence it is positive. By Remark 7,7 because $p_1 q_1' + q_1 q_1' > p_1 q_1' + q_1 p_1$ we deduce that x(.)x'(+)y(.)y'>x(.)y'(+)x'(.)y. Next we prove the property 7. Let x' as before but also satisfying x'∈w(x) that is x'<x. Then by property 5 we deduce that w(x)+y ⊆ w(x+y). Conversely let z<x+y z ∈w(x+y). Let the Cantor normal form of z be

$$z = \sum_1^n \omega^{\xi_i}(.)r_i$$

(we rearrange appropriately the normal forms of x, x', y, y', Z) with $r_i \in$ No. From the last inequality we get that in the lexiographical ordering it holds that $(r_i,...,r_n)<(p_1+q_1,...,p_n+q_n)$(*)

Let $k_i = \max_{1\le i\le n}\{r_i,q_i\}$ and $\lambda_i = \max_{1\le i\le n}\{r_i,p_i\}$

Let $z_1' = \sum_i^n \omega^{\xi_i}k_i$ and $z_2' = \sum_i^n \omega^{\xi_i}\lambda_i$ Then the following ordinals are defined: z₁'(-)z, z₁'(-)y, z₂'(-)z, z₂'(-)y, and also by Remark 7,7. It holds that z≤z₁' z≤z₂', y≤z₁' x≤z₂'. From the inequality (*) and the inequality (**) $q_i \le p_i + q_i$ i = 1,...,n and the definition of $k_i$ we infer that it holds in the lexicographical ordering, the inequality $(k_1,...,k_n) \le (p_1+q_1,...,p_n+q_n)$ similarly $(\lambda_i,...,\lambda_n) \le (p_1+q_1,...,p_n+q_n)$. Hence by Remark 7.7. it holds that z₁'≤x(+)y and z₂'≤x(+)y If for both z₁', z₂' holds that z₁'=x(+)y=z₂'.

Then $\overset{\max}{i}\{r_i, q_i\} = \overset{\max}{i}\{r_i, p_i\} = p_i + q_i$, hence $r_i = p_i + q_i$  $i=1,...,n$.

But then $z = x(+)y$ ,contradiction.

Let us suppose then, that $z_1' < x(+)y$.

Then if $z'' = z_1'(-)y$ by the last inequality we get that $z''(+)y = z_1'(-)y(+)y = z_1' < x(+)y$ or $z''(+)y < x(+)y$.

That is we proved that for every $z \in w(x+y)$ there is $z''$ an other ordinal with $z \leq z''(+)y < x+y$. If $z'' \geq x$ then $z''(+)y \geq x(+)y$ contradiction, hence $z'' < x$ that is $z'' \in w(x)$   and $z''(+)y \in w(x)+y$. From this and also that $w(x)+y \subseteq w(x+y)$ ,that we have already   proved ,we deduce that $w(x+y)$ and   $\{w(x)+y \cup x+w(y)\}$   are   cofinal sets;   we   write $cf(w(x+y)) = cf(\{w(x)+y \cup x+w(y)\})$ .  In  other words we haved proved the property 7.

Let us prove the property 8. As  we  have  already  remarked the difference is to be understood in the extension of  the  additive monoid $w(a)$ into the linearly ordered  Grothendick  group $k(w(a))$.  The  ordering  in $k(w(a))$  is  defined  by: $(x,y) \leq (x',y') \Leftrightarrow x+y' \leq x'+y$.

Where by $(x,y)$ we denote the equivalence class  of  the  free abelian group  generated  by  $w(a)$, which is denoted  by  $F_{a,b}(w(a))$ $(k(w(a)) = F_{a,b}(w(a))/([x+y]-[x]-[y]))$, in the process of taking the quotient by the normal subgroup generated by the elements of the form $[x+y]-[x]-[y]$ in $F_{a,b}(w(a))$ (the  corresponding generator of $x \in w(a)$, in $F_{a,b}(w(a))$ we denote by $[x]$), that  is defined by the representative $x+(-y)$. Needless to  mention that the natural difference in $w(a)$, isn't but an instance of difference in $k(w(a))$.

We   make   clear   that   $h^{-1}(\{h(x)(.)h(w(y))(+)h(w(x))(.)h(y) - h(w(x))(.)h(w(y))\} = \{v|v \in w(\alpha)$  and  $v = h^{-1}(h(x)(.)h(y)'(+)h(x')(.)h(y) - h(x')(.)h(y'))$  with  $x' \in w(x)$ $y' \in w(y)$ and $x,y \in w(a)\}$. By the property  6 we get that $h(x)(.)h(y)'(+)h(x')(.)h(y) < h(x)(.)h(y)(+)h(x')(.)h(y')$  hence

$h(x)(.)h(y)'(+)h(x')(.)h(y) - h(x')(.)h(y)' < h(x)(.)h(y)$

hence  $h^{-1}(\{h(x)(.)h(w(y))(+)h(w(x))(.)h(y) - h(w(x))(.)h(w(y))\}\} \subseteq w(x(.)y)$.

Conversely, let, $z \in w(x(.)y)$, that is $z < x(.)y$.

If $x(.)y = \sum_{1 \leq i,j \leq n} \omega^{\xi_i (+) \xi_j} (.)(p_i.q_i)$ then we also write for the normal form of z: $z = \sum_{1 \leq i,j \leq n} \omega^{\xi_i (+) \xi_j} (.)r_{ij}$ and $r_{ij} \in$ No. By Remark 7,7. We deduce that in the lexicographical ordering it holds that $(r_{11},...,r_{ij},...,r_{n,n}$ $) <$ $(p_1 p_1,...,p_i p_j,...,p_n . p_n)$. It is sufficient to prove that for every $(r_{11},...,r_{ij},...,r_{n,n}) < (p_1 p_1,...,p_i p_j,...,p_n.p_n)$ there are $(p_1',...,p_n')$ and $(q_1',...,q_n')$, $p_i, q_j \in$ No with $(p_1',...,p_n') < (p_1,...,p_n)$ and $(q_1',...,q_n') < (q_1,...,q_n)$ such that $(r_{11},...,r_{ij},...,r_{n,n}) \leq (p_1 q_1' + p_1' q_1 - p_1'.q_1',...,p_n.q_n' + p_n'.q_n - p_n'.q_n') < p_1.q_1,..., p_n.q_n)$.

But the property 8 holds for $a = \omega$, that is for the natural numbers. Hence there are $p_1'$, $q_1'$ with $p_1' < p_1$ $q_1' < q_1$ and $r_{11} \leq p_1 q_1' + p_1' q_1 - p_1'.q_1' < p_1' q_1$ and completing with arbitrary $p_i'$, $q_i'$ $i = 2,...,n$) that give positive the terms $p_i q_j' + p_i' q_j - p_i' q_j'$ (by elementary arithmetic of natural numbers this is always possible) we define

$x' = \sum_{1}^{n} \omega^{\xi_i} (.)p_i'$ and $y' = \sum_{1}^{n} \omega^{\xi_i} (.)q_i'$.

By the lexicographical ordering it holds that $x' \in w(x)$, $y' \in w(y)$ and $h(z) \leq h(x)(.)h(y')(+)h(x')(.)h(y) - h(x')(.)h(y') < h(x)(.)h(y)$. Hence the sets $W(x(.)y)$ and $h^{-1}(\{h(x)(.)h(w(y))$ $(+)$ $h(w(x))(.)h(y)$ $-$ $h(w(x))(.)h(w(y))\})$ are cofinal and we write

$cf(w(x(.)y)) = cfh^{-1}(\{h(x)(.)h(w(y))$ $(+)$ $h(w(x))(.)h(y)$ $-$ $h(w(x))(.)h(w(y))\})$.

This is the end of the proof of the property 8.    Q.E.D.

## Corollary 9 (first characterisation )

Every pair of operations in a principal initial segment of ordinal that satisfy the properties 0.1.2.3.4.5.6.7.8 of lemma 1 ,.is unique up-to-isomorphism and coincides with the Hessenberg natural operations .

**Remark**: The difference that appears in the property 8 is defined as in the remark after the proposition 8 .

**Proof**:Direct after the proposition 4 and 8      Q.E.D.

.

**Corollary 10.(Second characterisation )**

Every pair of field -inherited operations in a principal initial segment of ordinals, w(a)  that satisfy the properties  7. 8. coincides  with the natural sum and product of Hessenberg .

( For the existence of field-inherited operations in the ordinal   numbers see [C Conway J.H.1976 ] ch note pp 28 .)

**Proof**: The proof is immediate from  proposition 4 and 8.Q.E.D.

**Remark.11** It seems that N.L.Alling in his publications:

a)On the existence of real closed Fields that are $\eta_\alpha$ -sets of power $\omega_\alpha$ Transactions Amer.Math.Soc. 103 (1962) pp 341-352.

b)Conway's       field       of       surreal       numbers. Trans.Amer.Math.Soc.287 (1985) pp.365-386.

c)Fountations of Analysis over Surreal  number Fields .Math. studies 141 North-Holland 1987.

he is unaweare that if an initial segment of ordinals is contained in a set-field and it is cofinal with   the   field  ,(and  it  induces   the Hessenbeg operations   in it) then it has to be an initial segment of a principal ordinal that is  of type  $\omega^{\omega^\alpha}$ (see [Kutatowski K-Mostowski A.1968] ch VII §7 p. 252-253 exercises 1. 2.  3.)

Thus properties 0.1.2.3.4.5.6.7.8. can be taken as an axiomatic definition of the  Hessenberg  operations  without having to mention the non-commutative ordinal operations.

In a forthcoming paper, I  will  be  able  to prove      the      non-contradictory     of   properties 0.1.2.3.4.5.6.7.8. (actually the existence in Zermelo-Frankel set theory, of the operations +,.) without  using  the non-commutative  ordinal operations,neither field-inherited operations.  but through  transfinite  induction  and  other methods of universal algebra .

## List of special symbols

ω          : Small Greek letter omega,  the  first  infinit number.

α, b        : Small Greek letter alfa, an ordinal.

$\Omega_1$         : Capital  Greek  letter  omega  with   the subscript one.

$\aleph(x)$         : Aleph of x, the cardinality of the set X.

             N: the fisrt  capital  letter  of  the  Hebrew alfabet.

$\oplus, \circ$        : Natural sum and  product of G. Hessenberg plus and point
             in   parenthesis.

# Free algebrae and alternative definitions of the Hessenberg operations in the ordinal numbers . The ORDINAL NATURAL NUMBERS 2

## Constantine E. Kyritsis*

*Associate Prof. of the University of Ioannina, Greece.*
*ckiritsi@teiep.gr      C_kyrisis@yahoo.com, Dept. Accounting-Finance Psathaki Preveza 48100*

### Abstract.

It is proved and is given, in this paper, two alternative algebraic definitions of the Hessenberg natural numbers in the ordinal numbers: a) by definition with transfinite induction and two inductive rules , b) by the free algebras of the polynomial symbols of the commutative semirings with unit.

**Key Words:**Hessenberg natural operations**,** ordinal numbers,free algebras**,**semirings

**Subject Classification of AMS** *03,04,08,13,46*

**§0. Introduction** This is the second paper of a series of five papers that have as goal the definition of topological complete linearly ordered fields (continuous numbers) that include the real numbers and are obtained from the ordinal numbers in a method analogous to the way that Cauchy derived the real numbers from the natural numbers. We may call them linearly ordered Newton-Leibniz numbers. The author started and completed this research in the island of Samos during 1990-1992 .

We should not understand that with the current theory we suggest direct applications in the physical sciences. Not at all! **Matter is always finite**. Actually not even the real numbers are fully appropriate for the physical reality because they are based on the infinite too which does not exist in the material reality. This has been described in more detailed by the famous Nobel prize winner physicist E. Schrödinger in his book "Science and Humanity" (see [ Schrödinger E. 1961]. That is why the author has developed the digital or natural real numbers without the infinite with the corresponding Euclidean geometry and also Differential and Integral calculus, which is logically different from the classical. (See [Kyritsis 2017] and [Kyritsis 2019]) But the ordinal numbers and the surreal numbers reflect more the human **consciousness and perceptions** rather than properties of the physical material reality. Still such a discipline as the study of the continuum of the surreal numbers is an excellent spiritual, mental and metaphysical meditative practice probably better than many other metaphysical spiritual systems. ***It is certainly an active reminding to the scientists that the ontology of the universe is not only the finite matter but also the infinite perceptive consciousness.***

In this second paper on the same subject, I shall give two more ,and even simpler, algebraic characterizations of the Hessenberg natural operations in the ordinal numbers. These characterizations are actually alternative and direct definitions of the Hessenberg natural operations; independent from the standard non-commutative operations in the ordinal numbers .The main
results are the characterization theorems 4.7.

.These characterisations of the Hessenberg natural operations are :

a) *As operations defined by transfinite induction through two inductive rules that are already satisfied by the usual operations in the natural numbers .*

b) *An initial segment of a principal ordinal $\omega^{\omega^\alpha}$ in the Hessenberg natural operations is isomorphic with the free semiring of $\alpha$ many generators in the category of abelian semirings ;or isomorphic with the algebra of polynomial symbols of α indeterminates of the type of algebra of semirings with constants the natural numbers .*

The previous characterizations proves that the Hessenberg natural operations are the natural extensions in the ordinal numbers, of the usual operations in the natural numbers,. This turns out to be so, if we approach this subject from whatever aspect. Thus the Hessenberg natural operations should be coined as the standard abelian operations in the ordinal numbers, for all practical algebraic purposes .There are already extended applications of this . (see [ Conway J.H.])

The main application of the previous results is in the theory of ordinal real numbers (see [ Kyritsis C.E. ] ).The final result is that the three Hierarchies and different techniques of transfinite real numbers (see [ Gleyzal A. ]), of the surreal numbers (see [ Conway J.H.] ) ,of the ordinal real numbers (see [ Kyritsis C.E. ]) give by inductive limit or union the same class of numbers ,already known as the class No and to which we make reference in [ Kyritsis C.E. ] as the "***totally ordered Newton-Leibniz realm of numbers*** ".

# §1. The third algebraic characterizations of the Hessenberg natural operations.

Let a initial segment of an ordinal number of type $\beta=\omega^\alpha$. Let us define a binary operation ,denoted by $+$, in $W(\beta)$ by definition by transfinite induction (see e.g. [Kuratowski K.-Mostoeski A.] §4 pp 233, [Enderton B.H.], [ Frankel A.A.], [Kyritsis C.E.] Lemma 2, 3 ) and the inductive rule

$$p_+ : \cup \; w\!\left(\gamma\right)^{W\!\left(\omega^\lambda\right)x W\!\left(\omega^\gamma\right)} \rightarrow W\!\left(\gamma+1\right)$$
$$x,y < \alpha$$

defined by $p_+ (f)=S(\{f(W(x),y)\}\cup\{f(x,W(y))\})$ for an $f \in W\!\left(\gamma\right)^{W\!\left(\omega^\lambda\right)x W\!\left(\omega^\gamma\right)}$.

The definition by transfinite induction is supposed of two variables as is also the inductive rule (see [ Kyritsis C.E] Lemma 2,3). Thus, there exists a unique function denoted by $(+):W(\omega^\alpha)^2 \rightarrow W(\gamma+1)$, where the $\gamma$ is an ordinal number with

$\omega^\alpha \prec \gamma$ such that it satisfies the inductive rule $p_+$ ; thus it holds :

$p_+$       $x+y=S(\{x+W(y)\}\cup\{\{W(x)+y\})$ .

The restriction of this operation in $W(\omega)=\omega$ coincides with the usual operations of the natural numbers, since the addition of the natural numbers satisfies also the inductive rule $p_+$ :

**Lemma 0**. *The Hessenberg natural sum in a initial segment* $W(\omega^\alpha)$ , *satisfies the inductive rule* $p_+$ .

*Proof:* See [ Kyritsis C.E] Proposition 8; the arguments hold also for the initial segments of type $W(\omega^\alpha)$; if we are concerned only for the natural sum .                    Q.E.D.

Thus by the uniqueness ,the natural sum coincides with the operation defined with the inductive rule $p_+$.

**Corollary 1.** *The operation defined as before by the inductive rule* $p_+$, *satisfies the properties* 0.1.2.3.5.6.7.*(See [ Kyritsis C.E] lemma 1, the part of the properties that refer only to the sum ).*

*Proof:* See again the [ Kyritsis C.E] proposition 8.          Q.E.D.

Since the commutative monoid $W(\omega^\alpha)$ relative the Hessenberg natural sum satisfies the cancellation law, it has a monomorphic embedding in the Universal group of it, which at this case is called also the Grothendick group and it is denoted by $K(W(\omega^\alpha))$. Thus the difference $x-y$ is definable in $W(\omega^\alpha)$ with values in $K(W(\omega^\alpha))$.

See also [ Kyritsis C.E] the remark before the proof of the proposition 8.
Let an initial segment $W(\omega^\alpha)$ of a principal ordinal number $\omega^{\omega^\alpha}$ . Let the binary operation denoted by $(.):W(\omega^{\omega^\alpha})^2 \to W(\gamma)$ ,where the $\gamma$ is an ordinal number with $\omega^{\omega^\alpha}<\gamma$ defined with definition by transfinite induction and with inductive rule  the function

p

$$p: \bigcup_{x,y<\alpha} W(\gamma)^{W(\omega^\alpha) \times W(\omega^{\omega^\alpha})} \to W(\gamma+1)$$

such that for every

$$f \in W(\gamma)^{W(\omega^\alpha) \times W(\omega^{\omega^\alpha})}$$

p (f)=S({f(x,W(y))+f(W(x),y)-f(W(x),W(y))}∩ W(γ+1)) . Thus there is a unique function $(.):W(\omega^{\omega^\alpha})^2 \to W(\gamma+1)$ such that it satisfies the inductive rule p        x.y=S({W(x).y+x.W(y)-W(x).W(y)}∩W(γ+1)).

**Lemma   2.***The   Hessenberg   natural product  satisfies  the  inductive rule* p .

***Proof****:* See [ Kyritsis C.E]Proposition 8 .                    Q.E.D.

Therefore  by the uniqueness of the function (.) this operation coincides with the Hessenberg natural product .

***Corollary 3****. Let an initial segment of a principal ordinal number* $\omega^{\omega^\alpha}$ *. The operations that are defined  as before with the inductive rules* p₊, p *satisfy  the  properties* 0.1.2.3.4.5.6.7.8. *(See  [ Kyritsis  C.E] lemma 1)and  coincide with the Hessenberg natural sum and product .*

***Proof****:* See again [ Kyritsis C.E] proposition 8.

## Corollary 4. (third characterisation)

Let   an   initial   segment of a principal ordinal number  $\omega^{\omega^\alpha}$ . Two operations in W($\omega^{\omega^\alpha}$) are the Hessenberg natural operations  if  and only if they satisfy the inductive rules p₊, p.

Proof: Direct from lemma 2 and corollary 3.            Q.E.D.

## § 2 .The definition of the Hessenberg natural operations with finitary free algebras .

In this paragraph we shall prove a key result with respect to the Hessenberg  natural  operations.  We  shall  prove  that  the  Hessenberg operations are actually free finitary operations definable by the operations of the Natural numbers .

(see [ Graetzer G.] about operations of polynomial symbols ch 1 e.t.c.)

**Proposition 5.** Let an initial segment $W(\omega^\alpha)$ of an ordinal number of type $\omega^\alpha$. The commutative monoid $W(\omega^\alpha)$ relative to the Hessenberg natural sum is isomorphic with the commutative free monoid $\coprod_\alpha \mathbb{N}_0$ in the category of commutative monoids .

   **Remark**: The free monoid $\coprod_\alpha \mathbb{N}_0$ coincides with the polynomial algebra of polynomial symbols of the algebra of type $(\mathbb{N}_0,.)$,in other words of commutative monoids with nullary operations the constants of $\mathbb{N}_0$. Since the commutative monoids is an equational class (variety ) there are free commutative monoids ;(see [ Graetzer G.] ch 4 §25 corollary 2 pp 167 ).

   ***Proof***: Let us define a function h : $\coprod_\alpha \mathbb{N}_0 \to W(\omega^\alpha)$ by $h(x)=\omega^x$ for $x<\alpha$ and $h(n_1 x_1 +...+n_k x_k)=n_1 \omega_1{}^x +...+ n_k \omega_k{}^x$ for any $y=n_1 x_1 +...+ n_k x_k \in \coprod_\alpha \mathbb{N}_0$. The operations in the second part of the defining equation of h are the Hessenberg natural operations .By the definition of $\coprod_\alpha \mathbb{N}_0$ and the Cantor normal form of ordinal numbers in the Hessenberg operations we get that the h is 1-1 on-to and homomorphism of abelian monoids .Thus an isomorphism of commutative monoids . Q.E.D.

***Remark 6***. We deduce from the previous proposition that two initial segments $W(\omega^\alpha)$, $W(\omega^\beta)$ are algebraically isomorphic as commutative monoids if and only if $\aleph(\alpha)=$

$\aleph(\beta)$, in other words the ordinals $\alpha$, $\beta$ have the same cardinality .

**Proposition 7. (Fourth characterisation )** *Let an initial segment W($\omega^{\omega^{\omega^\alpha}}$ ) of a principal ordinal number $\omega^{\omega^{\omega^\alpha}}$ . The  commutative semiring W($\omega^{\omega^{\omega^\alpha}}$ ) relative to  the Hessenberg natural operations is isomorphic with the commutative free semiring*

*$N_0{}^{\left(\coprod_\alpha \mathbb{N}_0\right)}$ in the category of commutative semirings  with unit .*

**Remark** : The free commutative semiring with unit $N_0{}^{\left(\coprod_\alpha \mathbb{N}_0\right)}$ coincides with the polynomial algebra of polynomial symbols of the algebrae of type

$(N_0,+,.)$ in other words of the commutative semirings with nullary operations the constants of $N_0$. The commutative semirings with unit are an equational class thus they have free semirings;

(see again [ Graetzer G.] ch 4 § 25 corollary 2 pp 167 ).The semiring $N_0^{\left(\coprod_\alpha \mathbb{N}_0\right)}$ is constructed as the semigroup semiring of the semigroup $\left(\coprod_\alpha \mathbb{N}_0\right)$ written multiplicatively;(in a way analogous to the construction of the semigroup ring of a semigroup ).

**_Proof_**: Let us define a function as in the proof of proposition 5

$h_2 : N_0^{\left(\coprod_\alpha \mathbb{N}_0\right)} \to W(\omega^{\omega^\alpha})$ by $h_2(x)=\omega^{h(x)}$ for $x \in \coprod_\alpha \mathbb{N}_0$ where the h is as in the proof of the proposition 5 h : $\left(\coprod_\alpha \mathbb{N}_0\right) \to W(\omega^\alpha)$ and the $\left(\coprod_\alpha \mathbb{N}_0\right)$ is written multiplicatively; $y \in \mathbb{N}_0^{\left(\coprod_\alpha \mathbb{N}_0\right)}$ $y= n_1 x_1 +...+n_k x_k$ with $x_1,...,x_k \in \left(\coprod_\alpha \mathbb{N}_0\right)$ $h(y)=$

$n_1\omega^{h(x_1)} +...+ n_k\omega^{h(x_k)}$. Again by the definition of the $N_0^{\left(\coprod_\alpha \mathbb{N}_0\right)}$ and the uniqueness of the Cantor normal form in the Hessenberg natural operations (see [Kyritsis C.E.] Remark 7,5),b)) we get that the function h is an homomorphism of semirings,1-1 and on-to ;thus an isomorphism of abelian semigroups with unit . Q.E.D.

**Remark 8**.From the previous proposition and the dependence of the free semiring $N_0^{\left(\coprod_\alpha \mathbb{N}_0\right)}$, up-to-isomorphism,on the cardinality of the set α, we deduce that two initial segments $W(\omega^{\omega^\alpha})$ , $W(\omega^{\omega^\beta})$ are algebraically isomorphic relative to the Hessenberg natural operations if and only if $\aleph(\alpha)=\aleph(\beta)$; n other words if the ordinal numbers $\alpha$, $\beta$ are of the same cardinality .

# ORDINAL REAL NUMBERS 1. The ordinal characteristic.

## Constantine E. Kyritsis*

*\* Associate Prof. of the University of Ioannina. [ckiritsi@uoi.gr](mailto:ckiritsi@uoi.gr) [C_kyrisis@yahoo.com](mailto:C_kyrisis@yahoo.com),, Dept. Accounting-Finance Psathaki Preveza 48100*

## Abstract

In this paper are introduced the ordinal integers ,the ordinal rational numbers ,the ordinal real numbers ,the ordinal p-adic numbers ,the ordinal complex numbers and the ordinal quaternion numbers .It is also introduced the ordinal characteristic of linearly ordered fields. The final result of this series of papers shall be that the three different techniques of surreal numbers, of transfinite real numbers ,of ordinal real numbers give by inductive limit or union the same class of numbers known already as the class No and that would deserve the name the "infinitary totally ordered Newton-Leibniz realm of numbers ".

§0 **Introduction.** This is the third paper of a series of five papers that have as goal the definition of topological complete linearly ordered fields (continuous numbers) that include the real numbers and are obtained from the ordinal numbers in a method analogous to the way that Cauchy derived the real numbers from the natural numbers. We may call them linearly ordered Newton-Leibniz numbers. The author initiated and completed this research in the island of Samos in Greece during 1990-1992 .

We should not understand that with the current theory we suggest direct applications in the physical sciences. Not at all! **Matter is always finite**. Actually not even the real numbers are fully appropriate for the physical reality because they are based on the infinite too which does not exist in the material reality. This has been described in more detailed by the famous Nobel prize winner physicist E. Schrödinger in his book "Science and Humanity" (see [ Schrödinger E. 1961]. That is why the author has developed the digital or natural real numbers without the infinite with the corresponding Euclidean geometry and also Differential and Integral calculus, which is logically different from the classical. (See [Kyritsis 2017] and [Kyritsis 2019]) But the ordinal numbers and the surreal numbers reflect more the human **consciousness and perceptions** rather than properties of the physical material reality. Still such a discipline as the study of the continuum of the surreal numbers is an excellent spiritual, mental and metaphysical meditative practice probably better than many other metaphysical spiritual systems. ***It is certainly an active reminding to the scientists that the ontology of the universe is not only the finite matter but also the infinite perceptive consciousness.***

In a communication (1992) that the author had with N.L. Alling and his group of researchers on analysis on surreal numbers, suggested the term ordinal real numbers instead of surreal numbers. Some years later and before the present work appears for publication, it appeared in the bibliography conferences about *real ordinal numbers* .

In these last three papers is studied a special Hierarchy of transcendental over the real numbers, linearly ordered fields that are characterized by the property that they are fundamentally (Cauchy ) complete. It shall turn out that they are isomorphic to the transfinite real numbers (see [Glayzal A. (1937)]).The author was not familiar with the 5 pages paper of [ Glayzal A. (1937)] ,and his original term was "transfinite real numbers". When one year later (1991) he discovered the paper by A. Glayzal ,he changed the term to the next closest :"Ordinal Real Numbers" .One more year later he proved that the transfinite real numbers ,the surreal numbers and the ordinal real numbers were three different techniques leading to isomorphic field of numbers. He then suggested (1992) to researchers of surreal numbers, like N.L.Alling to use the more casual term "ordinal or transfinite real numbers " for the surreal numbers. In the present work it is introduced a new, better, classifying and more natural technique in order to define them. This technique I call *"free operations-fundamental completion"*. It is actually the same ideas that lead to the process of construction of the real numbers from the natural numbers through fundamental (Cauchy) sequenses. In the modern conceptual context of the theory of categories this may demand at least three adjunctions (see[ MacLane S 1971 ]).It is developed their elementary theory which belongs to algebra. Their definition uses the Hessenberg operations of the ordinal numbers .It may be considered as making use of an infinite dimensional K-theory which is mainly not created yet. In this first paper it is also introduced the ordinal characteristic of any linearly ordered field .It is a principal ordinal number, that is of type $\omega^{\omega^{\omega^{\lambda}}}$ . These numbers ,as defined with the present technique of the "free operations-fundamental completion " and prior to the proof that the resulting linearly ordered fields are isomorphic to the transfinite real numbers (as in [Glayzal A. (1937)]) ,we shall call Ordinal real numbers. The relevancy with the surreal numbers and the non-standard (hyper) real numbers ,shall be studied in a later paper. In detail, the next Hierarchies are defined:

1) The **Ordinal natural numbers**, denoted by $\mathbf{N_\alpha}$ .2) The **Ordinal integral numbers**, denoted by $\mathbf{Z_\alpha}$ 3)The **Ordinal rational numbers**, denoted by $\mathbf{Q_\alpha}$ 4)The **Ordinal p-adic numbers**, denoted by $\mathbf{Q_{\alpha,p}}$ 5)The **Ordinal real numbers**, denoted by $\mathbf{R_\alpha}$ 6)The **Ordinal comlpex numbers**, denoted by $\mathbf{C_\alpha}$7)The **Ordinal quaternion numbers**, denoted by $\mathbf{H_\alpha}$, of characteristic $\alpha$ . The fields $\mathbf{Q_{\alpha,p}}$, $\mathbf{R_\alpha}$, $\mathbf{C_\alpha}$, $\mathbf{H_\alpha}$ are fundamentaly (Cauchy)complete topological fields.

The field $\mathbf{R_\alpha}$ is also the unique maximal field of characteristic $\alpha$ ( that is, it is Hilbert complete) , and the unique fundamentally (Cauchy ) complete field of characteristic $\alpha$. It is also a real closed field , according to the theory of Artin-Schreier . These will be proved in the next paper on ordinal real numbers.

As it is known there are three more techniques and Hierarchies of transcendental over the real numbers, linearly ordered fields. Namely (in the historical order): The transfinite real numbers (see [Glayzal A. 1937 ]), and the surreal numbers (see [Conway J.H (1976) ]).

In this series of papers, it is proved (among other results ) that all the previous three different techniques and Hierarchies give by inductive limit, or by union, the same class of numbers (already known as the class **No** ).

## § 1. The ordinal characteristic of linearly ordered fields.

**Definition 0.** We remind the reader that a linearly (totally) *ordered, double abelian semigroup (semiring )* M is a set with two operations denoted by +,., such that with each one of them it is an abelian semigroup. Furthermore the distribution law holds for multiplication over addition. A linear ordering is supposed defined in M that satisfies the following compatibility conditions with the two operations 1) if x>y, x'>y' x,x',y,y'∈M then x+x'>y+y' and xx'+yy'>xy'+yx' (The symbol < is used for ≤ and not equal) if M is also a monoid relative to the two operations, and zero is absorbent unit for M, M is called *ordered double abelian monoid*. (semiring)

(e.g.The set of natural numbers ,denoted by N).

In the next we shall consider linearly (totally) ordered fields.( For a definition see [Lang S.] ch xi §1 pp 391).

Also in the next we shall use ordinal numbers. (For a reference to standard symbolism and definitions see [Kuratowski K.-Mostowski A. 1968] ch vii, [Cohn P.M. 1965] pp 1-36 )

In the following paragraphs we will not avoid the use of larger totalities than the sets of the Zermelo-Frankel set theory, namely classes.

We may suppose that we work in the Zermelo-Frankel set theory, augmented with axioms for classes also, as is presented for instance in bibliography [ Cohn P.M. 1965] p.1-36 with axioms A1-A11. Wee denote by $\Omega_1$ the class of the ordinal numbers. (The last capital letter of the Greek alphabet with subscript 1). The axioms A1-A11 allow for larger entities than sets, to define algebraic fields or integral domains or semi-groups. Hence we will also study classes that have two algebraic operations (Their Cartesian square treated as classes of sets of the form {{x, y}, {x}}, that is of ordered pairs) that satisfy the axioms of an algebraic field and have a subclass called the class of positive elements, with properties 1. 2., that they define a compatible ordering in the field (again as a class of ordered pairs) such classes that are ordered fields we will call again ordered fields and if we want to discriminate them from set-fields, especially when they are classes that are not sets, we will write for them that they are c-fields similarly we write c-integral domains or c-semigroups. We must not confuse the term "c-field" with the term "class-field" of the ordinary set-fields of "class-fields theory" (see [ Van der Waerden B.L 1970], [Artin E.-Tate J. 1967]). A subset (or subclass) denoted by X ⊆ F of a linearly ordered field F, is said to be <u>cofinal</u> with F, if for every a ∈ F there is a b ∈X with a≤b.

**Definition 4**. The ordinal characteristic is essentially a measurement of the size of a linearly ordered commutative field with a semi-ring of Ordinal natural numbers (Hessenberg natural commutative operations in the ordinal numbers, as developed in the two previous papers-sections). We embed systems of ordinal natural numbers in a linearly ordered field, so that not "gaps" exist. There is always a minimal such system the natural numbers themselves. The definition of the ordinal characteristic of such ordinal natural numbers (See Definition 6 below) is always the supremum of the ordinals which are contained in it, and it is a principal ordinal numbers as we have described in the previous paper-section. Then we embed with monomorphisms and with 1-1 functions , such semi-rings of ordinal natural numbers in a linearly ordered commutative field so that the 0 and 1 of the ordinal real numbers goes to the 0, 1 of the linearly ordered field and there are no "gaps", in other words the image is the minimal such possible set in the linearly ordered field. All such possible monomorphic with no gaps in a linearly ordered field, which is a set, give a set of corresponding ordinal characteristics of such semi-rings of ordinal natural numbers which is upper bounded, because of the cardinal and corresponding ordinal of the set and linearly ordered field. Thus as such ordinal are a subset of a well ordered set of ordinals it holds the supremum property, and there is such a supremum ordinal. Since also such a maximal embedding is also a semi-ring of ordinal natural numbers, this supremum is also a principal ordinal number which exist and its unique, it measures the size of the linearly ordered field and we call the it its ordinal characteristic. We say that the field (or integral domain or double abelian monoid ) F is of characteristic $\alpha$ and we shall write charF = $\alpha$.

If F is a c-field we include the case of characteristic $\Omega_1$ and we write charF = $\Omega_1$ if all ordinals contained in F is the class $\Omega_1$ and also it is a cofinal subclass with F.

**Remark** . In the case of a set-field F with $\alpha$ = charF, we do not need to suppose that the subset of elements of F corresponding to the ordinal in $\alpha$ by the definition 1 (it always exists ,by making use of the definition by transfinite induction and its version that uses only a set of functions sufficient for an inductive rule), see appendix A), is cofinal with F, as this is a consequence of the definition. For, if there is an element $X_0 \in F$ with $\beta < X_0$ for every ordinal number $\beta$ with $\beta \leq a$, then the set $\alpha \cup \{X_0\}$ can be extended , with the field operations ,to its closure in the natural Hessenberg operations (a semiring) (see [Kyritsis C. Alt] ) and it becomes similar to an initial segment of a principal ordinal number Thus $\alpha + 1$ is an ordinal contained in F, contradiction with the definition of a .

By the previous definitions we realize that every linearly ordered set-field has characteristic which is a limit ordinal number.

The fact that the linearly ordered field F has characteristic $\omega$ (the least infinite ordinal) is equivalent with the statement that the field F is Archimedean.

In the followings when we will work on a linearly ordered field denoted by F of ordinal characteristic $\alpha$, $\alpha$=charF (or $\Omega_1$= charF) we will supposed that is fixed an embedding of the ordinal numbers of the initial segment w($\alpha$) in the set F (or of $\Omega_1$ in F).

If the characteristic is $\omega$, the embedding is obviously unique as it can be proved by finite induction.

**Remark.5** Let a linearly ordered field denoted by F .Obviously there is an extension which is a real field .Let us denote by R(F) the real closure of F .(For results of the theory of Artin-Schreir on real and real closed fields see e.g.[ Lang S. 1984] ch xi .or [Artin E.-Shreier O. 1927]) Since R(F) can be obtained by adjunction of the square roots of the positive elements of F and Zorn's Lemma on algebraic extensions see[Lang S. 1984] ch i proposition 2.10 theorem 2.11 pp 397), it is direct that the characteristic of the real closure R(F) is the same with that of F.

For the definitions of the terms <u>infinite</u>, <u>finite</u>, <u>infinitesimal</u> elements in an extension of such fields, see e.g.[ Lang S] ch xi paragraph 1 pp 391, the definitions can be given relative to extensions of any linearly ordered field to an other linearly ordered field ,and not only extensions of the real numbers.

## §2 The ordinal natural numbers N . The ordinal- integers Z .

Let w($\alpha$) a principal initial segment of ordinal numbers. Let us denote by + and . the Hessenberg's natural sum and product in w($\alpha$). They satisfy properties 0.1.2.3.4.5.6. after lemma 1 in §1 in [ Kyritsis C.1991 Alter]

**Definition 6**. *The set w($\alpha$)=$\alpha$ where $\alpha = \omega^{\omega^x}$ for some ordinal x, is an abelian double monoid relative to sum and product and furthermore it satisfies the cancellation low (see [ Kyritsis C. 1991 Alter] lemma 1 ).This set I call the (double abelian) monoid of <u>ordinal natural numbers</u> of characteristic a and I denote it by $N_\alpha$. Thus $N_\alpha$ =$\alpha$.*

**Remark 7**. It is obvious that the (double abelian, well ordered ) monoid $N_\alpha$, is <u>the</u> minimal such monoid of characteristic $\alpha$ and the embedding of the ordinal numbers of W($\alpha$) in it is unique . Furthermore it can be proved by transfinite induction that it is a <u>unique factorization monoid</u> (called simply factorial monoid also).

The additive cancellation low in $\alpha$ has as a consequence that $\alpha$ is monomorphicaly embedded in its Grothendieck group denoted by k(α) (see [Lang   S.   1984] Ch.1 §9   p.   44). Furthermore the Grochendieck group $k(N_\alpha)$ can be ordered  by defining the set of positive elements $k(\alpha)^+ = \{v/v = (x,y)$ with $x,y \in w(\alpha)$ and $x > y\}$. We remind the reader  that if we  denote  by $F_{ab}(\alpha)$ the free abelian group generated    by    α,    and by   ((x+y)-x-y)   the   normal subgroup  of  $F_{ab}(\alpha)$ generated by  elements  of  the  form  (x+y)-x-y,

$$k(\alpha) \cong \frac{F_{ab}(\alpha)}{((x+y)-x-y)}$$

then

By (x,y) we denote the equivalence class that is defined in $F_{ab}(\alpha)$ in the   process  of  taking  the  quotient  group  $F_{ab}(\alpha)/((x+y)-x-y)$  by  the representative x+(-y).

The  first  part  of  property  6.  (lemma  1  in  [Kyritsis C.1991 Alter]) guarantees that this ordering in k(α) restricted on  α coincides with the usual ordering of ordinal numbers.

**Definition   8**. *The   ordered   Grothendieck   group  k(α)  of  an  initial segment of ordinals relative to natural sum, we call <u>transfinite cyclic group</u> of exponent α and we  denote it by $\Gamma_\alpha$. (by [Kuratowski K. Mostowski A. 1968] ch vii §7 pp 252-253 exercises 1.2.3.the ordinal α has to be of the type $\omega^x$. If the ordinal α is principal then I denote it also by $Z_\alpha$).*

Every element  of  the  group  $Z_\alpha$  is  represented  as  a difference x-y with $x,y \in w(\alpha)$. Then we define <u>multiplication</u> in $Z_\alpha$ by  the  rule

(*)  (x-y).(x'-y')=(x.x'+y.y')-(xy'+x'y)

where sum and product are  the  natural  sum  and  product  in w(α). This makes $Z_\alpha$ a commutative ring with unit (the element 1).

If (x-y)(x'-y') = 0 and both (x-y), (x'-y') are not zero, we get by property 6 in lemma 1 in [Kyritsis C. 1991Alter] that xx'+yy' ≠ xy'+yx' or (x-y)(x'-y') ≠ 0, contradiction. Then one of

(x-y), (x'-y') is zero that is  the  ring  $Z_\alpha$ has  no  divisors  of  zero and  it  is  an  integral  domain.  Remembering  that $Z_\alpha^+ = \{v|v \in Z\alpha$ and $v = (x,+y)$ with $x,y \in w(\alpha)$  $x > y\}$, by property 6 lemma 1 in [ Kyritsis C. 1991 Alter], we get that the sum and product of elements of $Z_\alpha^+$ are again elements  of  $Z_\alpha^+$.  From all these we get:

**Lemma   9**. *The  ring  $Z_\alpha$ is  a  linearly  ordered  integral  domain  of characteristic  the  principal  ordinal  α (see  §  1  Def.1).The set $Z_\alpha^+$ is a linearly ordered double abelian monoid and $Z_\alpha^+ \neq N_\alpha$*

**Definition      10** . *The      integral  domain  $Z_\alpha$  I  call  <u>ordinal  integers</u> of characteric α .*

The integral domain $Z_\alpha$ of characteristic $\alpha$ has <u>minimality</u> relative to its property of being an integral domain of characteristic α, in the following sense: Every integral domain of characteristic α contains a monomorphic image of $Z_\alpha$.

**Theorem 11** <u>(Minimality)</u>.

*Every integral domain $Z_\alpha$ is minimal integral domain of characteristic α. That is every integral domain of characteristic α, contains a monomorphic image of $Z_\alpha$.*

*Proof.* Put $R_\alpha$ an integral domain of characteristic α, where α is a principal ordinal number ($\alpha = \omega^{\omega^\lambda}$).

Then the initial segment w(α) is contained in $R_\alpha$ (more precisely an order preserving image of w(α)). The principal initial segment is closed to the integral domain operations and by theorem 13,14 of [ Kyritsis C. 1991 Alter], they coincide with the natural sum and product of Hessenberg. Then, applying the construction of this paragraph for the integral-domain $Z_\alpha$, we remain inside the integral-domain $R_\alpha$, that is $Z_\alpha \subseteq R_\alpha$. This proves the minimality.

**Remark 12**. The ordinal integers are semigroup-rings of quotient monoids of semigroups that are used to define as semigroup-rings the hierarchy of integral domains of the <u>transfinite integers</u> (see [Gleyzal A. 1937] pp 586).I use the term hierarchy not only as a well ordered sequence but also as a net (thus partially ordered ). The transfinite real numbers are thus an hierarchy.

The transfinite integers over the order-type λ symbolised by Z(λ), is the semigroup-ring (also module Z-algebra and integral domain) of the linearly ordered monoid $\sum_\lambda \mathbb{N}$ , where $\sum_\lambda \mathbb{N}$ is the coproduct, or direct sum denoted also by $\coprod_\lambda \mathbb{N}$ , of a family of isomorphic copies of N with set of indices the order-type λ. Thus Z(λ) =Z[ $\coprod_\lambda \mathbb{N}$ ]. Thus any ring of polynomials of a linearly ordered set of variables with integer coefficients is an integral domain of transfinite integers and conversely. It can be proved with the axiom of choice and transfinite induction , as in the case of finite set of variables, that Z(λ) is a <u>unique factorization domain</u> . On the other hand the Cantor normal form in the Hessenberg operations of the ordinal numbers (see lemma 6 in [Kyritsis C. 1991 Alter]) gives that any element x of $Z_\alpha$ is of the form $x = \omega^{x_1} y_1 + \ldots + \omega^{x_n} y_n, y_i \in Z, i = 1,\ldots,n, n \in \mathbb{N}$ $x_i$ are ordinals with $x_1 > \ldots > x_n$. The ordinal powers of $\omega$ in $Z_\alpha$ is an abelian well ordered monoid (see e.g. [Neumann B.H. 1949] §2 pp 204-205) of ordinal characteristic $\beta = \omega^x$, if $\alpha = \omega^{\omega^x}$. Let us denote it by $M_\beta$. Actually $M_\beta = \beta$. Let us denote by $\log_\circ(\alpha)$, or simply by $\lambda_\alpha$ the order type of the Archimedean equivalent classes of $M_\beta$. Then we get by the Cantor normal form that $Z_\alpha = Z[M_\beta]$ (The semigroup ring of $M_\beta$). The monoid $M_\beta$ can be obtained as quotient monoid of the free abelian multiplicative monoid of $\lambda_\alpha$ variables, which is the monoid $\coprod_{\lambda_\alpha} \mathbb{N}$ .

But $Z[\coprod_{\lambda_\alpha} \mathbb{N}] = Z(\lambda_\alpha)$, which was the assertion to be proved.

**Remark 13** The equation $Z_{\omega^{\omega^x}} = Z[M_{\omega^x}]$ gives an alternative, simpler definition of the ordinal integers without the use of the Hessenberg multiplication, since the ordinal powers of $\omega$ coincide n the abelian Hessenberg operations and the usual ordinal operations (see [Kyritsis C.1991 Alter] Remark 7.5) ) and without the use of the Grothentick group .The monoid $M_x$ is defined as the initial segment $W(\omega^x)$ (or simply as the ordinal $\omega^x$) in the Hessenberg addition .

## §3  The definition of the fields $Q_\alpha, R_\alpha, C_\alpha, H_\alpha$.

In this paragraph, I shall introduce the hierarchies of fields of ordinal rational ,real, complex ,quaternion numbers. These hierarchies give the unification of the other three techniques and hierarchies, namely of the transfinite real numbers, of the surreal numbers. Furthermore we introduced the hierarchies of transfinite complex and transfinite quaternion numbers.

**Definition 14.** *The localization (field of quotients) of the integral domain $Z_\alpha$, I will denote by $Q_\alpha$ and I will call <u>ordinal rational numbers</u> (of characteristic $\alpha$)  (see [Lang S. 1984] ChII §3).*

**Remark**. Since we have that cancellation low holds, we do not have to use the Malcev-Neuman theorem (see [Cohn P.M. 1965] Ch VII §3. Theorem 3.8). We define as set of positive element of $Q_\alpha$ the set $Q_\alpha^{+'} = \left\{ x \middle| x = \frac{m}{1} \text{ with } m,1 \in Z_\alpha \text{ and } m1 \in Z_\alpha^+ \right\}$. It is elementary in algebra that if the integral domain is linearly ordered then also its field of quotients (localization) with the previous definition for its set of positive elements, is a linearly ordered field with the restriction of its ordering on the integral domain to coincide with the ordering of the integral domain. Obviously the ordinals of the initial segment of $w(\alpha)$ are contained in $Z_\alpha$ and also in $Q_\alpha$. By a direct argument, holds also that the characteristic of $Q_\alpha$ is a: $\underline{\text{Char } Q_\alpha = \alpha}$.

**Remark** From the construction of $Q_\alpha$ we infer easily that $\aleph(Q_\alpha) = \aleph(\alpha)$ and if $\alpha < \beta$ where $\alpha$, $\beta$ are two principal ordinals then $Q_\alpha \subseteq Q_\beta$. The converse obviously holds.

**Lemma 15.** Every element x of the field $Q_\alpha$ is of the form $x = \dfrac{\omega^{\alpha_1} . a_1 + ... + \omega^{\alpha_n} . a_n}{\omega^{\beta_1} . b_1 + ... + \omega^{\beta_m} . a_m}$ where $\alpha_i$, $\beta_j \in w(\alpha)$ and $\alpha_1 > \alpha_2 > ... > \alpha_n \geq 0$, $\beta_1 > \beta_2 > ... > \beta_m \geq 0$ and $a_i$, $b_j$ for $i = 1,...,n$, $j = 1,...,m$ are finite integers.

*Proof.* Direct from the definition of localization and lemma 6 in [ Kyritsis C. 1991Alter].

**Theorem 17**. (Minimality)

*The field $Q_\alpha$ is a $\underline{\text{minimal}}$ field of characteristic $\alpha$, in the sense that every field of characteristic $\alpha$, contains the field $Q_\alpha$ (more precisely an order preserving monomorphic image of $Q_\alpha$).*

**Remark.** This property is already obvious for the field of rational numbers, that in the statement of Theorem 17 is denoted by $Q_\omega$.

*Proof.* Let a field of characteristic $\alpha$, that we denote by $F_\alpha$. Then the principal initial segment $w(\alpha)$ of ordinals is contained in $F_\alpha$ and the field-inherited operations coincide with the natural sum and product of Hessenberg (see theorem 14 in [ Kyritsis C. 1991 Alter]). Then constructing first the integral domain $Z_\alpha$ and afterwards its localization $Q_\alpha$ we always remain in the field $F_\alpha$.

Thus $Q_\alpha \subseteq F_\alpha$ (or more precisely $h(Q_\alpha) \subseteq F_\alpha$ where h is a order-preserving monomorphism of $Q_\alpha$ in to $F_\alpha$) Q.E.D.

**Definition 18**. *The (strong) Cauchy completion of the topological field $Q_\alpha$ we denote by $R_\alpha$ and I call $\underline{\text{ordinal real numbers}}$ of characteristic $\alpha$.*

**Remark.**The process of extensions ,beginning with a principal initial ordinal $\alpha = N_\alpha$ which is $\underline{\text{the}}$ minimal double, abelian monoid of characteristic $\alpha$, and ending with the field $R_\alpha$ which is $\underline{\text{the}}$ maximal field of characteristic $\alpha$ ,we call $\underline{\text{K-fundamental densification}}$ .

**Lemma 19**. *The characteristic of the (strong) Cauchy completion of a linearly ordered field F ,is the same with that of the field F.*

*Proof.* If the characteristic of the field is $\alpha$, let us denote it by $F_\alpha$, and its completion by $\hat{F_\alpha}$. Obviously the characteristic of $\hat{F_\alpha}$ is not less than $\alpha$.

Suppose that there is an ordinal $\beta$ with $\alpha < \beta$ which is contained in $\hat{F_\alpha}$ (see Definition 1). Then there is a Cauchy net $\{x_i | i \in I\}$ of elements of $F_\alpha$ that converges to $\lim\limits_{i \in I} x_i = \beta$. Let $\varepsilon \in F_\alpha$ $0 < \varepsilon < 1$, then there is $i_0 \in I$ such that for every $i \in I$ $i \geq i_0$

$x_i \in (b - \varepsilon, b + \varepsilon)$. But this gives an element of $F_\alpha$ greater than $\alpha$, hence than every element of $F_\alpha$, which is a contradiction. Thus Char $R_\alpha = \alpha$. Q.E.D.

**Corollary 20**. *The characteristic of $R_\alpha$ is $\alpha$.*

From the definition of $R_\alpha$ we infer that $\aleph(R_\alpha) \leq 2^{\aleph(\alpha)}$ and that $\alpha < \beta \Leftrightarrow R_\alpha \overset{\subset}{\cdot} R_\beta$ for two principal ordinals denoted by $\alpha, \beta$.

**Remark.21** We denote by $R(\lambda)$ the transfinite real numbers of order-base $\lambda$. It holds by definition that $R(\lambda) = R((LR^\lambda))$, where $LR^\lambda$ is the lexicographic product of a family of isomorphic copies of the real numbers $R$, with set of indices the order-type $\lambda$.

**Remark**. It is said that a field $F$ has formal power series representation, if there is a formal power series ring $R((G))$ and a ideal $I$ of it such that $F$ has a monomorphic image in $R((G))/I$. From the universal embedding property of the hierarchy of transfinite real numbers we get that every linearly ordered field has formal power series representation .Thus:

**Corollary 22.** *The fields of ordinal real numbers R, have formal power series representation ,with real coefficients.*

**Definition 23** *The field $C_\alpha = R_\alpha[i]$ I call <u>ordinal complex numbers</u> of characteristic $\alpha$.*

**Definition 24**. *The field $C(\lambda) = R(\lambda)[i]$ we call <u>transfinite complex numbers</u> of base-order $\lambda$. Actually it is the field $C(\lambda) = C((LR^\lambda))$.*

**Definition 25**. *The quaternion extension field of the field $R_\alpha$ (or of $C_\alpha$) by the units i, j, k with $i^2 = j^2 = k^2 = ijk = -1$, I call the <u>ordinal quaternion numbers of characteristic $\alpha$</u> and I denote them by $H_\alpha$. They are non-commutative fields (following the terminology e.g. of A.Weil in [ Weil A. 1967]) that are transcendental extension of the non-commutative field H of quaternion numbers.*

**Definition 26**. *The formal power series fields $H(\lambda) = H((LR^\lambda))$ we call <u>transfinite quaternion</u> numbers of base-order $\lambda$.*

For a proof that H((LR$^\lambda$)) is a (non commutative ) field see [Neumann B.H.1949] part I.

## §4 The ordinal p-adic numbers Q$_{\alpha,p}$ .

As it is known  if F is a linearly ordered field ,and K a linearly ordered subfield of the real numbers and F$\mid$K is an extension respecting the ordering, then this extension defines the order-valuation (see [N.L.Alling 1987] ch 6 § 6.00 pp 207) .Actually every extension of any two linearly ordered fields  F, K, K$\subseteq$F, respecting the ordering, defines a place, thus a valuation v. (I use the place and valuation as  are defined e.g. by O.Zariski in [Zariski O.-. Samuel P.1958] vol ii ch vi §2, §8.and not as are defined by A.Weil in [Weil A. 1967] ch iii or by v.der Waerden in [Van der Waerden B.L. 1970] vol ii ch 18 .The definition of Zariski is equivalent with the definition of v.der Waerden only for the non Archimedean valuations of the latter).

The place-ring is the F$_v$ ={x/x $\in$ F and there are a, b $\in$K with a<x<b }. The maximal ideal of the place (or valuation v ) is the ideal of infinitesimals of K relative to F.

This valuation we call underline{extension - valuation} (and the corresponding place underline{extension - place}) It  has as special case the order valuation .The rank of the extension- place (see [Zariski O.-. Samuel P.1958] vol.II §3 pp 9) we call the underline{rank of the extension} .If char(F)>char(K) then the extension is transcendental ,and has transcendental degree and basis ;the latter is to be found in the ideal of infinitesimals or in the set of infinite elements .

**Definition 27 .**

*Let F a field of ordinal characteristic. Let R a subring of F that has F as its field of quotients. Let p  a prime ideal of R, such that the triple (pR$_p$, R$_p$, F) where R$_p$ is the localization of R at p, defines a place of F.  Such  a  place  (or valuation denoted by v$_p$) I call underline{p-adic} of the field F. In  the  valuation  topology  of  the  valuation  v$_p$, that has a local base of zero the ideals of R ) the field F is a topological field and the (strong) Cauchy completion I denote by F$_p$, it is a (topological field ) and I  call  underline{p-adic extension field of F}.*

**Definition  28**.*For   F=Q$_a$ and   R=Z$_a$ in   the   previous   definition the   field  Q$_{a,p}$ I  call  underline{ ordinal p-adic numbers}  of characteristic α.*

**Final remark** .Using inductive limit ,or union of  the elements of the hierarchies of the previous ordinal and transfinite number systems, we get corresponding classes of numbers .The classes of ordinal natural, integer, rational, real, complex, quaternion numbers denoted respectively by Ω$_1$, (or On ), Ω$_1$Z, Ω$_1$Q, Ω$_1$R, Ω$_1$C, Ω$_1$H.

And the classes of transfinite integer, rational, real, complex, quaternion numbers denoted respectively by:

CZ, CQ, CR, CC, CH.

**Acknowledgments.** I would like to thank professors W.A.J.Laxemburg and A.Kechris (Mathematics Department of the CALTECH) for the interest they showed and that they gave to me the opportunity to lecture about the ordinal real numbers in CALTECH. Also the professors H.Enderton and G.Moschovakis (Mathematics Department of the UCLA) for their interest and encouragement to continue this project.

List of special symbols

$\omega$          : Small Greek letter omega, the first infinit number.

$\alpha, \beta$     : Small Greek letter alfa, an ordinal.

$\alpha = \omega^{\omega^{\lambda}}$     Ordinal alpha $\alpha$ equal to omega in the power of omega in the power of x

$\Omega^0$       : Capital Greek letter omega with the superscript zero.

$F^a$       : Capital letter with superscript a. The of algebraic elements of a field F.

char F    : The characteristic of a field denoted by F.

$\cong$        : Equiralence relation of Commensurateness.

$\sim$        : Equiralence relation of comparability.

tr.d.(x) : The transcendance degree initial of words tr.(anscendance) and d.(egree).

N(x)    : Aleph of x, the cardinality of the set X. N: the fisrt capital letter of the Hebrew alfabet.

cf(X)=cf(Y) :  The sets x and Y are cofinal.

W($\alpha$)    : Initial segment of ordinal naumbers defined by the ordinal number a.

$\oplus$,$\circ$    :  Natural sum and product of G. Hessenberg plus and point in parenthesis.

N$\alpha$,Z$\alpha$,Q$\alpha$,R$\alpha$,: Double-lined capital letters with subscript small Greek letters

C$\alpha$,H$\alpha$        namely transfinite positive integers, intergals, rationals reats, complex and quatenion numbers.

$Z\alpha_1{}^*\omega$    : The dual lually compact abelian groups of the transfinite integers Za. The capital letter Z double-lined wiuth subscripts two Greek let-$\alpha$ (alpha) and $\omega$ (omega) and superscript a star

$T_\alpha$    : Transfinite circle groups: Capital letter T with subscript a small Greek letter.

$^*$X, $^*$R et.c : A non-standard enlergement structure capital letter X with left superscript a star.

$\xi$No    : A sureal number field of characteristic $\xi$. A small Greek letter followed by the symbol No.

C,RC$^*$R,No : The            c-structures (classes) previous symbols following the capital

CN,CZ,CQ,.    latin letter C

CC,CH

$\hat{X}$ : Strong Canchy competition of a topological space capital letter with cap.

$\Sigma$ : Capital Greek letter sigma symbol for summation.

$\overset{\bullet}{D}_\alpha$ : The open full-linary tree of leight a. Capital latin D with subscript a small Greek letter and in upper place a small zero.

The ordinal real numbers 1. The ordinal characteristic.

APPENTIX A.

**A MORE EFFECTIVE FORM OF DEFINITION BY TRANSFINITE INDUCTION.**

1.Given a set Z and an ordinal $\alpha$, let $\Phi$ be a set of $\xi$-sequences with the properties:

a) If f belongs to $\Phi$ then $f/W(\xi)$ belongs to $\Phi$ for every $\xi <=$ domain of f.

b) For every $\xi<\alpha$ there is at least one f belonging to $\Phi$ with $\xi=w(\xi)=$domain(f) and values belonging to Z.

c)If $f_\xi$ is an $\alpha$-sequence of $\xi$-sequences of $\Phi$ such that whenever $\gamma<\xi_1$

, $\xi_2 <\alpha$ , $f_{\xi 1}/w(\gamma) = f_{\xi 2}/w(\gamma)$ ;then the $\alpha$-sequence $c_\alpha(\xi)=f_\xi(\xi)$, belongs to $\Phi$ also.

For each function h in $Z^\Phi$, there is one and only one transfinite sequence f defined on $\xi<=\alpha$,

f in $\Phi$ and such that $f(\xi)=h[f/w(\xi)]$ for every $\xi<=\alpha$ .

The function h is called **a recursive rule for $\Phi$.** The set $\Phi$ with the properties a). b), c), is called ,sufficient for recursive rules.

*Proof*: Not much different than the ordinary form of definition by transfinite induction.

# ORDINAL REAL NUMBERS 2. The "Cartesian" arithmetization of order types.

## Constantine E. Kyritsis*

* Associate Prof. of University of Ioannina, Greece. *ckiritsi@uoi.gr* *C_kyrisis@yahoo.com*, Dept. Accounting-Finance Psathaki Preveza 48100

### Abstract

In this paper the main results are :Proofs that the ordinal real numbers are real closed fields and complete up-to-characteristic .They are also Dedekind ,and Archemedean complete fields .They are real formal power series fields and Pythagorean complete fields It is proved and discussed the K-fundamental arithmetisation and the binary arithmetisation of the order types .

**Key words:**Real closed commutative fields,Grothendick group,Archemidean complete fields,linearly ordered commutative fields,full binary trees

**Subject Classification of AMS** *03,04,08,13,46*

§0 **Introduction** .

The author initiated and completed this research in the island of Samos in Greece during 1990-1992.

We should not understand that with the current theory we suggest direct applications in the physical sciences. Not at all! **Matter is always finite**. Actually not even the real numbers are fully appropriate for the physical reality because they are based on the infinite too which does not exist in the material reality. This has been described in more detailed by the famous Nobel prize winner physicist E. Schrödinger in his book "Science and Humanity" (see [ Schrödinger E. 1961]. That is why the author has developed the digital or natural real numbers without the infinite with the corresponding Euclidean geometry and also Differential and Integral calculus, which is logically different from the classical. (See [Kyritsis 2017] and [Kyritsis 2019]) But the ordinal numbers and the surreal numbers reflect more the human **consciousness and perceptions** rather than properties of the physical material reality. Still such a discipline as the study of the continuum of the surreal numbers is an excellent spiritual, mental and metaphysical meditative practice probably better than many other metaphysical spiritual systems. *It is certainly an active reminding to the scientists that the ontology of the universe is not only the finite matter but also the infinite perceptive consciousness.*

In this second paper on ordinal real numbers are proved, the main (elementary) properties of them. It is proved that the ordinal real numbers $R_\alpha$ of characteristic $\alpha$, is the maximal field of characteristic $\alpha$ (maximality) and that it is , according to the theory of Artin-Screier, a real closed field. (It turned out ,after the work was completed and by thinking aside, that they are also Archimedean complete (see [ Glayzal A. 1937]),formal power series fields with real coefficients ,Dedekind complete (see [Massaza,Carla 1971]), and Pythagorean fields ).

It is also proved a classification theorem which is analogous to the Hölder theorem for the Archimedean linearly ordered fields.In particular it is proved that any linearly ordered field of
characteristic $\alpha$ contains the field $Q_\alpha$ of ordinal
rational numbers of characteristic $\alpha$, as a dense subfield and it is contained in the field $R_\alpha$ of ordinal real numbers of characteristic $\alpha$, as a subfield ..As it is known, the linear segments of elementary euclidean geometry can be defined as special order-types with Archimedean property, and Archimedean (Hilbert) completness through axioms (see e.g. for a not ancient approach the Hilbert axiomatisation in [ Hilbert D.1977] ch 1 ).It can be proved to be order isomorphic with subsets of the real number field R. This is well-known and it can be called, the elementary arithmetisation of the order-types of Euclidean linear segments. On this fact is based the Cartesian idea of analytic geometry. This was an important turning point in the developments of the ideas and techniques of mathematics, of the discrete nature of numbers and continuous nature of geometry. The basic principle is that *the continuum is developed from the discrete* and not vice versa! An instance of this principle is the development of images and animation in computers through pixels and bits! It is surprising that in one of the consequences of the theory of ordinal real numbers, it is proved a far more advanced and complete result for the whole category of order types that has as corollary the previous important and elementary arithmetisation. Although more advanced, the result remains in the context of elementary theory of ordinal real numbers .In this result any order type can be "discretised" or "arithmetised" through the ordinal numbers.

The process of definition of the maximal fields $R_\alpha$, from the minimal (double well ordered ) monoids $N_\alpha = \alpha$, of principal ordinal numbers, we call K- fundamental densification. It is proved that any order-type is order isomorphic to a subset of some field $R_\alpha$. Thus any order-type is constuctible by K-fundamental densification from ordinal numbers .This is called the K-arithmetisation of the order-types. Although in the way it is presented, this result is softly obtained, throws new light to the relation of ordinal numbers and order-types ,this relation turns out to be similar to the elementary relation of numbers and line segments in geometry. Also it, holds a second kind of arithmetisation ,the binary arithmeti sation which we state in the same paragraph .

## §1 On the topology of linearly ordered fields. Local deepness, α-sequences.

The ordering of any linearly ordered field F defines a well known topology : the order-topology denoted by $T_<$. In this topology, as it is known, the field F is a topological field.

This topology has very good separation properties; it is a $T_1$-$T_5$ topology, that is a completely normal topology (see for instance [ Lynn A.Steen-Seebach J.A. Jr 1970] § 39 p. 66-68, also see [ Munkress J.R. 1975] Chapters I, II)

The previously described order topology is also called the locally convex topology compatible with the order (see [Nachbin L. 1976] Ch I, II). (The convexity defined by the order).

**Definition 1.** Let X be a topological space. Let p ∈ X. The least ordinal α such that it exists a (local) base denoted by Bp of open neighbourhoods of the point p which is an α-sequence such that if x<y<α, $U_x \subset U_y$, is called local deepness of X at p.

We notice that the concept of local deepness is very close to the concept of local weight of a topological space, where instead of ordinal we have an initial ordinal that is a cardinal number (see [Kuratowski K. 1966] V-I p. 53-54).

Examples of topological spaces such that every point has local deepness, are the $\xi^*$-uniform topological spaces as they are defined in [Cohen L.W.Goffman C. 1949] pp 66 conditions 1.2.3.4.

As in the case of fields that are classes, we may permit topological spaces that are classes and the open sets is a class of subclasses closed to union and finite intersection. For such spaces, the local deepness may be $\Omega_1$ that is the class of all ordinal numbers.

**Proposition 2.** Let X a topological space and α, a limit ordinal such that every point has local-deepness α let A⊆X. It holds that $x \in \bar{A} \Leftrightarrow$ there is an β-sequence {$x_s$|s<β} from elements of A such that $\lim_{s<\beta} x_s = x$. In other words topological convergence in X can be treated with β-sequences where β=car(α) is the upper character of α (see [N.L.Alling 1987] ch 1 §1.30 pp 29)

The proof is almost direct and to save space we shall not give it here.

**Proposition 3.** Let a field denoted by F of ordinal characteristic α, where α is a limit ordinal. Then every point x ∈F in the order-topology has local-deepness car(α), where car(α) is the upper character of α (see [ N.L.Alling 1987] ch 1 §1.30 pp 29 ).

The proof is again direct and outside the scope of the paper.

**Corollary 4.** Convergence in the order-topology of a field of ordinal characteristic α, can be treated with β-sequences β≥car(α)

Needless to say, that in the case in which the topological space is a class and the local deepness is $\Omega_1$, then convergence can be treated with $\Omega_1$-sequences.

### §2 The Holder- type classification .

**Lemma 5.** In every field of characteristic αthe field **Qα** is a dense subset.

*Proof.* Let a field of characteristic α, which we denote by Fα. By the theorem 17 of [Kyritsis C. OR1] the field Qα is a subfield of Fα. Let us suppose that it is not dense in Fα. Then there are two elements x,y ∈ Fα x<y , such that there is no element of Qα in the internal [x,y]. Then the element z = y-x is Qα-infinitesimal.

This holds because $Q\alpha = \overline{L}(x) \cup \overline{R}(x) = \overline{L}(y) \cup \overline{R}(y)$ where $\overline{L}(x) = \left\{ v | v \in Q\alpha \ \ v \leq x \right\} \overline{R}(x) = \left\{ v | v \in Q\alpha \ \ x \leq v \right\}$ and $\overline{L}(y)$ $\overline{R}(y)$ similarly. But by the hypothesis $\overline{L}(x) = \overline{L}(y), \overline{R}(x) = \overline{R}(y)$. *Thus* $Q\ \alpha = \overline{L}(x) \cup \overline{R}(y)$ and every element of Qα can be written as $r_2$-$r_1$ where $r_2 \in \overline{R}(y)$ and $r_1 \in \overline{L}(x)$. Also we have that $0 < y-x < r_2-r_1$. Then y-x is a Qα-infinitesimal and $\frac{1}{y-x}$ is a Qα-infinite element of Fα, thus $\frac{1}{y-x} > \alpha$, contradiction since Char Fα = α.

Then there are not two element y,x ∈ Fα x<y with no element of Qα in [x,y], and Qα is dense in Fα. Q E D

**Remark**. Thus every field Fα of characteristic α is a <u>Weil completion</u> of the field Qα of ordinal rational numbers (see [ Weil A.] ChIII Definition 2 but applied not only to local fields).

**Theorem 6** (<u>Maximality</u> or <u>completness up-to-characteristic</u> ).

The field Rα is <u>the maximal</u> field, of characteristic α. In the sense that every field of characteristic α is contained as subfield of Rα (more precisely Rαcontains an order preserving monomorphic image of the field).The field $R_\alpha$ is the <u>unique</u> fundamentally complete field of characteristic α.

**Remark**. This theorem is analogous to the well-known Holders theorem theorem ςηιψη στατες that every linearly ordered Archimedean field is a subfield of the field of real numbers. In other words the field of real numbers is the maximal Archimedean linearly ordered field. The previous property of the ordinal real numbers $R_\alpha$ relative to their characteristic ,we call <u>maximality</u> or <u>completeness up-to characteristic</u> .

But as an erroneous application of terms R is also the minimal Cauchy complete field of characteristic ω and this also applies for the fields Rα in the sense that a completion of a linearly ordered field of characteristic α must be the field Rα .

*Proof.* Let any field of ordinal characteristic α denoted by Fα. By theorem 17 of [ Kyritsis C. 1991], the field Qα is contained in Fα: $Q\alpha \subseteq F\alpha$. Let x $\in$ Fα. Let (L(x), R(x)) be the cut that x defines on Qα (L(x) = {v|v $\in$ Qα v<x}, R(x) = {v|v $\in$ Qα x<v}). Since Qα is dense in Fα (Lemma **5**). There is a Cauchy α-sequence {$x_n$|n $\in$ w(α)} of elements of Qα that converges in Fα to x (all topologies are the order-topologies). Hence $Q\alpha \subseteq F\alpha \subseteq R\alpha$ and the field Rα is a maximal field of characteristic α ; but also the field Rα is actually a minimal Cauchy complete field of characteristic α in the sense that the (strong) Cauchy completion $\hat{F\alpha}$ of any field $F_\alpha$ of characteristic α contains the field $R_\alpha$:$Q_\alpha \subseteq F_\alpha$ has as a concequence that $R_\alpha \subseteq \hat{F\alpha}$ . Thus if $F_\alpha$ is complete then $R_\alpha \subseteq F_\alpha$, $F_\alpha \subseteq R_\alpha$ hence $F_\alpha = R_\alpha$ Q.E.D.

The theory of Artin-Schreier of real closed fields has an excellant application to the ordinal real numbers .

**Corollary 7.** The fields of ordinal real numbers Rα are real closed fields.

Proof. Direct from Theorem 6 , and remark 5 of [ Kyritsis C.1991] Q.E.D.

**Post written Remark A.** The author developed the theory of ordinal real numbers during 1990-1992 He had used the name "transfinite real numbers" without being aware that this term had been introduced by A.Glayzal during 1937 for his theory of linearly ordered fields. From the moment he fell upon the work of A.Glayzal (see [ Glayzal A. 1937 ]) in the bibliography of the Book of N.L alling (see [N.L.Alling 1987 ] ) he changed the title to "Ordinal Real Numbers" . After the work had been completed ,the author realised , by thinking aside, a quite unexpected and not unhappy fact :That the fields of ordinal real numbers are algebraically and order isomorphic to the fields of transfinite real numbers of Galyzal .This can be deduced by the fact that the fields of transfinite real numbers are exactly all the Archemidean complete fields (see [Glayzal A. 1937] theorems 4,8,9) and by the maximality of the ordinal real numbers (theorem 6). Thus if $R_\alpha$ is a field of ordinal real numbers of characterisic $\alpha$, any Archemidean (linearly ordered field ) extension of it ,it shall have the same characteristic with $R_\alpha$. It seems that it can be proved , that any cofinal (coterminal) linearly ordered field extension ,is of the same characteristic . By the maximality of $R_\alpha$ (theorem 6) it shall have to coinside with $R_\alpha$. In other words the fields of ordinal real numbers are Archemidean complete fields (although they may be non-Archemidean ).But this is a characteristic property of the fields of transfinite real numbersb of Glayzal.

Thus they are order and field isomorphic with fields of transfinite real numbers .Conversely ,let any field $R(\lambda)$ of transfinite real numbers of Archemidean base $\lambda$. Let us denote by $\alpha$ its ordinal characteristic .Let us suppose that there is an order and field extension of it with the same characteristic .Then it has to be an Archemidean extension of $R(\lambda)$. By the Archemidean completness of the transfinite real numbers ,it has to coinside with the $R(\lambda)$. Thus the transfinite real numbers are also complete up-to-characteristic .

But this is a characteristic property of the fields of the ordinal real numbers.Hence they are order and field isomorphic with fields of ordinal real numbers .Thus the ordinal real numbers should be considered as a different technique ,nevertheless indispensable and more far reaching .It is the technique that everyone would like to work.

**Post written Remark B** .Let a field $R_\alpha$ of ordinal real numbers of ordinal characteristic $\alpha$. It is also a field of transfinite real numbers of archemidean base $\lambda$. The set of all elements of $R_\alpha$ that as formal power series have support of ordinality less than $\beta \leq o(\lambda)$=maximum ordinality of well ordered set of $\lambda$, and which we denote by $R_{\alpha,\beta}$ is a field ,subfield of $R_\alpha$ .Indeed $R_{\alpha,o(\lambda)} = R_\alpha$. For the applications and especially with measurment proceeses ,the fields $R_{\alpha,\omega}$ are of prime interest and indispensable .

**Post written remark C** .The facts of the previous remark ,have as a concequence that the fields of ordinal real numbers are formal power series fields with coefficients in the real numbers and exponents in some order types.Thus the n-roots of their positive elements are contained in them (see [Neumann B.H. 1949] pp 211 ,4.91 Corollary).In other words they are Pythagorean complete fields.

**Theorem 8.** (The Holder-type classification theorem).

Every field of ordinal characteristic α, denoted by Fα (where α is a principal ordinal) is contained between the fields Qα and Rα : $Q_\alpha \subseteq F_\alpha \subseteq R_\alpha$ .

*Proof.* Contained in the proofs of the theorem 7 and lemma 5    Q.E.D.

**Remark.9** The previous theorem gives that the hierarchy of ordinal real numbers has *universal embedding property* for the category of linearly ordered fields, that is every linearly ordered field has an monomorphic image in some field of the hierarchy.The hierarchy of transfinite real numbers is known to have, also, this property .Such hierarchies we call universal embedding hierarchies. Especially the hierarchy of ordinal real numbers after the classification theorem 8 ,we call also, universal classification hierarchy.

**Remark.10** We notice that since every order type λ is order-embeddable in some transfinite real number field R(λ) (see [Glayzal A. 1937] )as Archemidean base which in its turn is embeddable in some ordinal real number field $R_\alpha$ ,the above two hierarchies as hierarhies of order-types are universal embedding hierarchies for the category of order-types .Let an order type λ ; the least principal ordinal number α such that λ is order-embeddable (by a monomorphism) in the order-type and field $R_\alpha$, is called the fundamental density of the order type λ and is denoted by df(λ).

**Remark.** In the [ Massaza Carla, 1971] Definition I , is defined which cuts are the Dedekind cuts in linearly ordered fields .It is proved also that the Dedekind completion D(F) of a linaerly ordered field F is also its Cauchy completion (in the order topology ).If we take the Dedekind completion $D(R_\alpha)$ of a field of ordinal real numbers $R_\alpha$, it has to be its Cauchy completion which is again the $R_\alpha$. Thus the fields of the ordinal real numbers are also Dedekind complete . Conversely ,let any Dedekind complete linearly field F .Let us denote with α its ordinal characteristic .Then by the Holder type classification (theorem 8 ) it is a subfield of the field $R_\alpha$ of ordinal real numbers of characteristic α .Since the Dedekind completion D(F) =F coincides with the F and also with its Cauchy completion ,we get that F=$R_\alpha$, because the Cauchy completion of F is the $R_\alpha$. In other words the class of Dedekind complete fields coincides with the class of the fields of ordinal real numbers .

Summarising we mention that the fields of ordinal real numbers have at least four kinds of completnesses that characterise them : Cauchy completeness ,Dedekind completeness,completeness up-to-characteristic, Archemidean completeness .It seems that he previous four completnesses can be summarised by saying that there is no cofinal (coterminal ) order field extensions of them ;in short they are _cofinally complete_ ,or _cofinally maximal_ .They are also realcomplete (closed ,Artin-Shreier ) and Pythagorean complete.

Remark. By corollary 7 we get that the field $C_\alpha$ is the algebraic closure of $R_\alpha$ : $C_\alpha = R_\infty^{\lambda}$ .

We close this paragraph by mentioning that an axiomatic definition of the field Rα (α is a principal ordinal) would be the following:

### First axiomatic definition of Rα.

*The field of ordinal real numbers Rα is the unique Funtamental (Canchy)-complete, in the order-topology, field of characteristic α.*

### Second axiomatic definition of $R_\alpha$ .

*The field of ordinal real numbers $R_\alpha$ is the unique complete (up-to-characteristic) field of characteristic α . These definitions apply even in the case of the field of real numbers (a = ω).*

### § 3 The arithmetisation of order-types .

**Remark**.As it is known the linear segments of elementary Euclidean geometry can be defined as special order-types with Archimedean property and Hilbert completness through axioms (see e.g. for a not ancient approach the Hilbert axiomatisation in [Hilbert D 1977] ch 1 ).Then ,they can be proved to be order isomorfic with subsets of the real number field R. This is known as the elementary arithmetisation of the order-types of Euclidean linear segments.

### Proposition 10.(the K- funtamental arithmetisation theorem of order-types.)

*Every order-type λ is K-arithmetisable with ordinal numbers and has a fundamental density df(λ) which is a principal ordinal number .*

In the next paper ,after the unification theorem of the transfinite real ,surreal ,ordinal real numbers ,a second arithmetisation theorem shall be proved. Two more universal hierarchies of formal power series fields shall be, also, proved that they are universal embedding hierarchies .We state these results here. For the definition of tree ,height of a tree, level of a tree, binary tree e.t.c.see [Kuratowski K.-Mostowski A. 1968] ch ii § 1, § 2 . The binary tree of height the ordinal α we denote with $D_\alpha$. After the previously mentioned unification theorem 17 of the next paper we get that the hierarchy of binary trees is a universal embedding hierarchy for the order- types . Since the binary trees are subsets of linearly ordered fields and their elements consisting exclusively from 1's in the binary sequence, correspond to the ordinal numbers with the Hessenberg operations (see also [Conway J.H. 1976] ch 3 note pp 28 and also [ Kyritsis C. 1991Alt] the characterisation theorem ) this universal embedding                    property we call also underline{binary arithmetisation} .The least ordinal α such that an order-type λ is order embeddable in the binary tree $D_\alpha$ ,we call the binary density of λ ,and we denote it by  db(λ).

**Theorem 11 ( The binary arithmetisation theorem of order-types )**

Every order-type λ  is binary arithmetisable and has a binary density db(λ) which is an ordinal number .

From the previous theorem ,by denoting a level of height α of a binary tree ,by $T_\alpha$ ,and giving to the Cartesian product $\prod_{\beta<\alpha} T_\beta$  the lexicographical ordering ,we also get the next  :

**Corollary 12**. The formal power series hierarchies $R((D_\alpha))$, $R_{\beta<\alpha}\left(\left(\prod T_\beta\right)\right)$, are universal embedding hierarchies for the linearly ordered  fields .

**§ 4 Some general results on linearly ordered fields .**

In this paragraph we give some results generally for the category of linearly ordered fields. To save space we shall not give the proofs, since they do not have serious dificulties,nevertheless we shall indicate how they can be obtained .

**Lemma 13** (On the rank and characteristic)

*Let us suppose that the characteristic of the field F is  $\omega^{\omega^\alpha}$  where α, is α limit ordinal. It holds  that  the  rank  of the extension F/K is a cofinal order-type  with  the  characteristic  of  the  field  F.  That is  cf(r(F/K))= cf(charF)=cf(char F-char K).*

Remark.For the definition of the rank of an extension see [ Kyritsis C. 1991] § 4. For the proof of the previous theorem we use the existence  for  any  principal  ordinal  $\omega^{\omega^\alpha}$  of  the  ordinal  real  numbers fields  $R_{\circ^{\omega^\alpha}}$ of characteristic  $\omega^{\omega^\alpha}$ .

Let F be  a linearly ordered field. If x ∈F by L(x) we denote the set L(x) = {y| y ∈F y<x} and by R(x) the set R(x) ={y| y ∈F x<y}.

By elementary arguments on linearly ordered  fields  the following identities can be proved.

**Lemma  14**

Let x, y ∈ F. The following hold

1.  L(-x) = - L(x)                    R(-x) = -R(x)

2.  $L(x+y) = L(x)+y = x+L(y)$

$R(x+y) = R(x)+y = x+R(y)$

3.  $L(x.y) = L(x).y + xL(y) - L(x).L(y) = R(x).y + xR(y) - R(x)R(y)$

$R(x.y) = L(x).y+xR(y) - L(x).R(y) = R(x).y + xL(y) - R(x).L(y)$

4.  $\left.\begin{array}{c} y < R(x) \\ L(y) < x \end{array}\right\} \Leftrightarrow y \leq x$

5.  $L\left(x^{-1}\right) = \dfrac{1+(R(x)-x)L(y)}{R(x)} = \dfrac{1+(L(x)-X)R(y)}{L(x)}$

$R\left(x^{-1}\right) = \dfrac{1+(L(x)-x)L(y)}{L(x)} = \dfrac{1+(R(x)-x)R(y)}{R(x)}$

The previous identities show also that the definition of operations used to define the surreal number fields are not something peculiar to these fields but hold in any <u>linearly ordered field</u> .

In the next paper of this work we will understand the true peculiarity of the technique of the surreal numbers.

**Lemma 15** If F/k is an extension of two linearly ordered fields , it holds that

tr.d.(F/k) $\leq 2^{\aleph(\text{Char.F})}$ where tr.d.is the transandental degree of the extension .

**Remark**. For the definition of the transandental degree, base e.t.c see for instance [ Zariski O.-Samuel P. 1958] vol. I pp. 95-102 also [Kyritsis C. 1991 ] § 4 ). The proof is obtained by using the Holder-type classification for F :$Q_\alpha \subseteq F \subseteq R_\alpha$ where α=char(F).

The next proposition shows that all the information of an extension of linearly ordered fields is to be found in the ideal of infinitesimals (or in the infinite elements). **Proposition 16**. Let F/k, F'/k two (ordered) extensions of the same linearly ordered field k. If the ideals of K-infinitesimals of the extension denoted respectively by $m_F$ and $m_{F'}$ are isomorphic as ordered integral domains ,then this isomorphism is extendable to an algebraic and order isomorphism of the fields F, F'.

**Remark** .The proof is direct from the definitions.

**Remark**. An extension F/k of the linearly ordered field k to F, is transcendental if Char F>char k and then the field F is an infinite dimensional vector space over k.

**Proposition 17** . Let F be a linearly ordered field of characteristic char(F)= $\omega^{\omega^\alpha}$ where α is a limit ordinal . It holds that the field F in the order topology is <u>totally disconnected</u> .

**Remark**.        The        proof uses the existence, for every principal ordinal $\omega^{\omega^\alpha}$, of the fields of ordinal real numbers R $\omega^{\omega^\alpha}$ .

**Theorem 18** The classes of transfinite real numbers CR, and of ordinal real numbers $\Omega_1 R$, coinside.

*Proof.* Since both Hierarchies of transfinite real and ordinal real numbers have the  universal embedding property (see remark 9 ) ,every transfinite real number-field is contained in some ordinal real number-field and every ordinal real-number field in some transfinite real        number-field.Thus CR $\subseteq \Omega_1 R$  and  $\Omega_1 R \subseteq$ CR, and CR = $\Omega_1 R$ .    Q.E.D.

## § 5    The A-Archimedeanity

The, at least  two different, definitions of archemideanity, that can be found for instance in [Glayzal A. 1937] and in other authors as in [ Conway J.H. 1976 ] or [ Arin E. Schreier O. 1927] give us the opportunity to treat them in unified way through the concept of archemideanity relative to a monoid.

The fact that the linearly ordered field  F  has characteristic ω (the  least  infinite  ordinal) is equivalent  with  the statement that the field F is Archimedean according to any (classical) known definition.

Let us denote by G a linearly ordered group and by A a monoid of endomorphisms of G  as a group.

It is said that x is <u>A-Archimedean</u> to y where x,y $\in$ G  iff there are a,b $\in$ A  with  a(x)≥y and b(y)≥x. If A is the domain Z of integers (the endomorphisms are  multiplication  with  an integer )we simply say that <u>x is Archimedean to y</u>. If for every pair x,y of elements of G holds that x is A-Archimedean to y, it is said that <u>G is A-Archimedean</u>

Let F be a linearly ordered field .If we consider it as an additive group, and we denote by $A_1$ a monoid of endomorphisms of the additive group , we get the concept of x being A-<u>additively Archimedean</u> to y. If we consider the multiplicative group $F^*$ and we take a monoid, denoted by $A_2$, of edomorphisms of the multiplicative group, we get the concept of x  being A-<u>multiplicatively Archimedean</u> to y.

Let A=$A_1 V A_2$ be the monoid of mappings from F to F generated by the previous monoids . It is said that x is A-<u>field-Archimedean</u> to y iff there are a, b $\in$ A such that a(x)≥y, b(y)≥x.

In any extension F/k of a field K by a field F, where F,k are fields of ordinal characteristic with char F>Char K, if we take as $A_1$, to be the multiplication with elements from the field K ( considering the field F as a linear space over K), we get the concept of x being K-additively Archimedean to y.( For K=R this is also known as "x is commensurate to y " see [Conway J.H. 1976] ch 3 pp 31 ).

If $A_1$ is the multiplication with integers and $A_2$ is power with integral exponents ,then it is simply said that x is field Archimedean to y (Known also from the A. Gleyzal's definition of Archimedeanity)

A non-Archimedean linearly ordered field denoted by F is simply a linearly ordered field for which not all pairs (x,y) of its elements are mutually additively Archimedean. (Thus charF>ω ) But it can be very well A-additively Archimedean for other monoids A.In particular if charF=α and A is the monoid of endomorphisms of the additive group of F defined by (field ) multiplication with ordinals less than α, then it is A-additively Archimedean and we denote it by writing that it is α-additively Archimedean

**Acknowledgments.** I would like to thank professors W.A.J.Luxemburg and A. Kechris (Mathematics Department of the CALTECH) for the interest they showed and that they gave to me the opportunity to lecture about the ordinal real numbers in CALTECH. Also the professors H. Enderton and G.Moschovakis (Mathematics Department of the UCLA) for their interest and encouragement to continue this project.

## List of special symbols

$\alpha,\beta,\omega$ : Small Greek letters

$\Omega_1$ : Capital Greek letter omega with the subscript 1

$F^a$ : Capital letter F with superscript a.

N : Capital Aleph ,the first letter of the hebrew alphabet . In the text is used a capital script. letter n .

$\oplus,\circ$ : cross in a circle, point in a circle .

$N\alpha, Z\alpha, Q\alpha, R\alpha$,: Roman capital letters with subscript small Greek letters

$C\alpha, H\alpha$

$^{*}X$, $^{*}R$ et.c : Capital standard or roman letters with left superscript a star.

CN,CZ,CQ, :Capital standard letter c followed by capital letters

$C^{*}R$, with possibly a left superscript a star

$\hat{X}$ : Capital tstandard letter with a cap.

$\Sigma$ : Capital Greek letter sigma

$\dot{D}_\alpha$ : Capital standard D with subscript a small Greek letter and in upper place a small zero.

# ORDINAL REAL  NUMBER 3. The continuum of the transfinite real , surreal, ordinal real, numbers ; unification .

## Constantine E. Kyritsis*

*Associate Prof. of the University of Ioannina.* ckiritsi@uoi.gr C_kyrisis@yahoo.com, *Dept. Accounting-Finance Psathaki Preveza 48100*

## Abstract

In this last paper on the theory of the ordinal real numbers, is proved ,that the three different techniques and hierarchies of transfinite real numbers , of the surreal numbers ,of the ordinal real numbers, give by inductive limit or union the same class of numbers.

**Key words:** Linearly ordered commutative fields,transfinite real numbers,surreal numbers,formal power series fields

**Subject Classification of AMS** *03,04,08,13,46*

§ **0 Introduction** . In this third paper on ordinal real numbers, *it is proved that the three different techniques and Hierarchies of transfinite real number-fields, of surreal numbers, and of ordinal real numbers , give by inductive limit ,or union, the same class of numbers and continuum, ,already known as the class No of surreal numbers*.It can be characterized, simply, as the smallest (linearly ordered field  which is a ) class and contains every linearly ordered set- field as a subfield . This class ,and also the category of linearly ordered set-fields, we call   **" the linearly ordered transfinite continuum  of infinite numbers".** It is obvious that without the set theory of G. Cantor as it is formalized, for instance, by Zermelo-Frankel and a correct thinking about the infinite, this "realm of  numbers" would   not be definable.

We should not understand that with the current theory we suggest direct applications in the physical sciences.  Not at all!  **Matter is always finite**. Actually not even the real numbers are fully appropriate for the physical reality because they are based on the infinite too which does not exist in the material reality. This has been described in more detailed by the famous

Nobel prize winner physicist E. Schrödinger in his book "Science and Humanity" (see [ Schrödinger E. 1961]. That is why the author has developed the digital or natural real numbers without the infinite with the corresponding Euclidean geometry and also Differential and Integral calculus, which is logically different from the classical. (See [Kyritsis 2017] and [Kyritsis 2019]) But the ordinal numbers and the surreal numbers reflect more the human **consciousness and perceptions** rather than properties of the physical material reality. Still such a discipline as the study of the continuum of the surreal numbers is an excellent spiritual, mental and metaphysical meditative practice probably better than many other metaphysical spiritual systems. ***It is certainly an active reminding to the scientists that the ontology of the universe is not only the finite matter but also the infinite perceptive consciousness.***

It is not directly apparent that so different techniques and ideas would have such an underlying unity. It is, also , surprising that, although the Hessenberg operations were very early known in the theory of ordinal numbers, (at least since 1906, see [ Gleyzal A. 1937]) no one went far enough to define through them, fields in a way similar to the way that the real numbers are defined from the natural numbers. Although G.Cantor, himself was conceiving the ordinals as a natural continuation of the natural numbers (see [Frankell A. A. 1953 ] introduction pp 3 ) ,as it is kmown, he rejected the attempts to define infinitesimals through them . (see [Frankell A. A. 1953] ch ii § 7.7 pp 120). We could speculate that un underlying reason for this, might be that, his set-theory was already strongly attacked and was facing the danger of final rejection ,and these were good enough reasons to avoid the additional charge that his theory " opened the door" to infinitesimals . In spite of this, there are many who might consider that although the present results are coming now, nevertheless it is too late and ,they might speculate ,for this long delay (more than eighty years) and diversion of ideas and technique, nevertheless on the same subject, we could suspect systematic obstructions, that came outside the mathematics. Nevertheless, there are others who consider that it is too early for such a development ,and especially for an analysis on such numbers. It seems that it has never been published any **"partially ordered transfinite continuum of infinite numbers"** (in other words a category of transcendental extensions of the real numbers ,that are partially ordered fields and complete in the order topology ) with reasonably "good" properties for a classification.

In this paper we use the surreal numbers, as they are definable in the Zermelo-Frankel set theory, through the binary trees, directly as a class, and not as union of some set fields.(The original technique of J.H.Conway).I met J.H.Conway during 1992 at Philadelphia in the USA, I talked to him about the new developments in this area of research and I gave to him the present work but as he told me he had more than a decade that for the last time he had active interest in the subject. I is somehow necessary to make use of classes instead of sets; since, for the kind of "induction" that the J.H.Conway uses, we prove that it is reduced to the usual transfinite induction on the height of the elements of the trees; but in their union as a class and <u>not</u> for each one of them separately as a set; in the latter case in which the trees are sets the induction fails .The key-point is to prove that for every cut that J.H.Conway uses it does really <u>exist</u> a <u>unique</u> element of the trees of least height . "simplest number" as it is used to be called ). This is a very crucial point, for the whole technique of the surreal numbers, to work, and it seems that it has been obscured, by not paying sufficient attention to it

The author has initially included also the non-standard real numbers in the classification. As they are also linearly ordered fields and the present classification is of all linearly ordered fields it was natural to include them. There were experts in non-standard analysis that were glad about it. Nevertheless there were experts that insisted that according to the initial definition of A.Robinson and not of later definitions, it was not claimed that the non-standard real numbers were sets inside Zermelo-Frankel system. Only if Zermelo-Frankel system was used to model meta-mathematics also the they would be also sets. This was nevertheless different as such sets would models of meta-mathematical entities different than the sets that are models of mathematical and not meta-mathematical entities. Because of their arguments and in spite the fact that this made some other researchers of non-standard mathematics unhappy, the author prefers in this first publication about ordinal real numbers not to include the non-standard real numbers in the unification. Any definition nevertheless that has the non-standard real numbers as ordinary sets of Zermelo-Frankel set theory, would naturally lead to a straightforward proof that such fields are always subfields of some field of ordinal real numbers! The author has already produced pages with this proof that is based on the premise that I mentioned.

## § 2. The surreal numbers .

In this paragraph we define the class No of surreal numbers inside the ZF-set theory.We use the binary trees (see [ Conway J.H. 1976] appendix to part zero pp 65 and [Kuratwski K.-Mostowski A.

1968] Ch IX §1, §2).The crucial point is to prove that for the cuts defined by J.H.Conway in these trees it does really <u>exist</u> a <u>unique</u> element strictly greater than all the elements of the left section and strictly smaller than all the elements of the right  section (the  "simplest  number" ).Through this the Conway-induction us reduced to the usual transfinite induction on the height of the elements of the tree .As we shall see this works for the union of all trees as a class but fails for each one set-tree .For the definition of the tree, binary tree, height, levels of the tree ,$H_\xi$-set see  [Kuratwski K.-Mostowski A. 1968]  Ch  IX  §1, §2 Theorem 2, . The binary tree of height α we denote by $D_\alpha$ . More precisely we are  interested for the trees of the next definition.

**Definition 1**. *Let α be an ordinal . We define* $\overset{\bullet}{D}_\alpha = \{x | x \in Da$ *such that there is β<α such that for the element x as a zero-one sequence x = $\{x_\xi | \xi < a\}$ holds that $x_\beta = 1$ and $x_\xi = 0$ for  ξ>β\}.*

We call the set $\overset{\bullet}{D}_\alpha$ <u>the open full-binary tree</u> of height α.

We also remind that if for the height α, holds that $\aleph(\alpha)$ is a cofinal to α regular aleph: $\aleph(\alpha) = \aleph_{cf(\alpha)} = \aleph_\xi$ the open full-binary tree is an $H_\xi$ set,  (see [Kuratwski K.-Mostowski A. 1968] ChIX §2 Theorem 2,  the proof works also  for  trees $\overset{\bullet}{D}_\alpha$ where $\aleph(\alpha) = \aleph_{cf(\alpha)}$

**Lemma 2**. *For every pair of subsets L, R of  the  open-full-binary tree* $\overset{\bullet}{D}_\alpha$ *of height the  ordinal  α, such that $\aleph(\alpha)$ is a regular aleph, and holds    that:    for    every    l∈L,    r∈R,    l    <    r,    and $\aleph(L), \aleph(R)$ < $\aleph(\alpha)$), there  is  exactly  one  element $x_0$ of least height in $\overset{\bullet}{D}_\alpha$ such that l < $x_0$ < r for every  l ∈ L, r ∈ R.*

*Proof.* Let D(L) = $\{x | x \in \overset{\bullet}{D}_\alpha$ such that there exists l ∈ L with x ≤ l} and I(R) = $\{x | x \in \overset{\bullet}{D}_\alpha$ such that there exists r ∈ R with r ∈ x} that is D(L), I(R) are  the  decreasing  and  increasing  lower  and  upper  half subsets of $\overset{\bullet}{D}_\alpha$ determined by L, R, in the  linear ordering  of $\overset{\bullet}{D}_\alpha$ as a tree (see [Kuratwski K.-Mostowski A. 1968] Ch IX §1  Lemma A). Let the set M = $\{x | x \ e \ \overset{\bullet}{D}_\alpha$ and for every $\not\subset$ 1 ∈ D(L), r ∈ I(R) it holds that l < x < r}. By the $H_\xi$ property of $\overset{\bullet}{D}_\alpha$  it  holds that M ≠ ∅ . Let A = {β|β is an ordinal number such that there is x ∈ M with x ∈ $T_\beta$ where  $T_\beta$ is the  β-level of Dα in other words there is x ∈ M of  height β}. Let $\alpha_0$ =  min  A. Let $D\alpha_0(L)$ $I\alpha_0(R)$ the subsets of D(L) R(L) of elements  of  height less than $\alpha_0$, and let $M\alpha_0 \subseteq M$ the subset of M that consists of elements of height $\alpha_0$. Suppose

that the set $M\alpha_0$ contains two elements x, y with e.g. $x \leq y$. We will prove that $M\alpha_0$ contains only one element.

Let $x'=\{x_\beta|\beta<\alpha_0\}$ that is that part of the $\alpha_0$-sequence x with terms of indifes less than $\alpha_0$. And the same also with $y' = \{y_\beta|\beta < \alpha_0\}$. Then there is $l_x$ or $r_x$ and $l_y$ or $r_y$ respectively in $D\alpha_0(L)$, $I\alpha_0(R)$ such that they are equal with x', y'. If $x=r_x$ then, if the $\alpha_0$-term of x is 0 or 1, in both cases $x > r_x$, contradiction. Hence there is no such $r_x$ and also such $r_y$. Then $l_x=x'$ $l_y=y'$ and $l_x \leq l_y$. The $\alpha_0$-term of x and y might be 0 or 1. The only possible cases are $\{x = (l_x,0), y = (l_y,0)\}$, $\{x = (l_x,0), y = (l_y,1\}$ $\{x = (l_x,1), (l_y,1)\}$, $\{x = (l_x,1), y = (l_y,0\}$ where with the parenthesis we symbolize the $\alpha_0$- sequence which is the elements x, y. Let us suppose that $x \neq y$ and $D\beta(x) = \left\{x_\beta|\beta < \delta\right\}$, the part of the $\alpha_0$-sequence with terms with indices less than $\delta$, with $\delta \leq \alpha$. Let the least value of $\delta$, be denoted by $\delta_0$ such that $D\beta(x) = D\beta(y)$, $\delta_0 \leq a_0$ and $X\delta_0 \neq Y\delta_0$. If holds that $0 = X\delta_0 \neq Y\delta_0 \neq 1$ because $x<y$. In the sequent, let $z=(D\beta(x)=D\beta(y)$ $\beta < \delta_0$, 1). Then $x < z \leq y$. If $\delta_0 = a_0$ then $x=y$ because $X\alpha_0 = Y\alpha_0 = 1$. Then $\delta_0 < \alpha_0$ and also $z < y$ and $x < z <$ y and the height of z is $\delta_0 < a_0$ contradiction. Hence $x=y$, and $M\alpha_0$ contains only one element. It also holds that if we restrict to $D_c(L)$, $I_c(R)$ where $D_c(L) = DC(L) \cap \overset{\bullet}{D}_c I_c(R) = I(R) \cap \overset{\bullet}{D}_c$ (and L, R have height $<c$), then $a_0 \leq c$ by the $H_\xi$-property of $\overset{\bullet}{D}_c$ if c is also such that $\aleph(\alpha) = \aleph_{cf(\alpha)}$ Q.E.D.

**Definition 3**. *The open full binary tree $\overset{\bullet}{D}\alpha$ of height $\alpha$, such that $\aleph(\alpha)$ is a cofinal to $\alpha$, regular aleph , I call  <u>regular open  full- binary tree</u>.*

The property of the previous lemma of a regular open full-binary tree I call <u>$H_\xi$-leveled Dedekind completness</u>.

We remark that the class of regular alephs is unbounded (see [Kuratwski K.-Mostowski A. 1968] p. 275 relation 5 ) Thus the class of ordinals $\alpha$ such that $\aleph(\alpha) = \aleph_{cf(\alpha)}$ is unbounded.

The next definition is the definition of the class of surreal numbers in the ZF-set theory and it depends as we mentioned on the lemma 2 .As it is seen ,in the hypotheses of the lemma 2 the cardinality of halfs of the cut is bounded by $\aleph(\alpha)$. If it is to include all possible cuts of the tree $\overset{\bullet}{D}\alpha$ then the lemma 2 will give the element $x_o$ in some tree $\overset{\bullet}{D}\beta$, of sufficient greater height,thus outside the original tree $D_\alpha$. This is why we mentioned that  the definition of surreal numbers (with the original technique of J.H.Conway ) does not apply to the trees $\overset{\bullet}{D}\alpha$ separately .

**Definition 4**. Let U$\overset{\bullet}{D}\alpha$=No be the union of all regular open full-binary trees. It is a class (after axiom A2.(see [Cohn P.M. 1965] p1-36)) Operations may be defined in this linearly ordered class according to the formulae of Lemma 2 in [Kyritsis C.1991 Alt. or Free etc.)] II, that hold for every linearly ordered field that is:

1. let α be an ordinal with $\aleph(\alpha)=\aleph_{cf(\alpha)}$ and L,R subsets of $\overset{\bullet}{D}\alpha$ such that for every $l \in L$, $r \in R$ holds that $l < r$. Then there exists a regular aleph β such that $L,R \subseteq \overset{\bullet}{D}_\beta$ and $\aleph(L)$, $\aleph(R)< \aleph(\beta)$. Then there is by lemma 2 a unique element $x_0 \in \overset{\bullet}{D}_\beta$ of least height such that $l < x_0 < r$ for every $l \in L$, $r \in R$, we denote this element by {L|R} and we write $x_0$= {L|R}. We note that although $L,R \subseteq \overset{\bullet}{D}_\alpha$, it holds that $x_0 \in \overset{\bullet}{D}_\beta$ and α <β.

2. If $x,y \in \overset{\bullet}{D}\alpha$ and we denote the height of x, y by h(x), h(y) and by L(x), L(y), R(x), R(y) the sets

$$L(x)=\left\{\nu \middle| \nu \in \overset{\bullet}{D}\alpha \quad h(\nu)<h(x) \text{ and } \nu < x\right\} L(y)=\left\{V \middle| V \in \overset{\bullet}{D}\alpha, h(\nu)<h(y) \text{ and } \nu < y\right\},$$

$$R(x)=\left\{\nu \middle| \nu \in \overset{\bullet}{D}\alpha, h(\nu)<h(x) \text{ and } x < \nu\right\} R(y)=\left\{\nu \middle| \nu \in \overset{\bullet}{D}\alpha, h(y)<h(x) \text{ and } y < \nu\right\}$$

Then the operations are defined through simultaneous two-variable transfinite induction in the form of the lemma 2,3 in [ Kyritsis C. 1991 Free etc.], for the heights of the trees $\overset{\bullet}{D}\alpha$ where for the initial segments of ordinals we substitute the corresponding trees of No (For every ordinal β<α such that $N(\beta)=N_{cf(\beta)}$ corresponds a tree $\overset{\bullet}{D}\beta$ ). Thus the function of operation is defined not on $w(\alpha)^2$ but on $\overset{\bullet}{D}\alpha^2$. For the addition, the next rule is used x+y={L(x)+y ∪ x+L(y)| x+R(y)} ∪ R(x)+y}.

3. The opposite is defined by:

-x = {-R(x)|-L(x)}

4. Multiplication is defined by

x.y={L(x).y+xL(y)-L(x).L(y)∪R(x).y+xR(y)-R(x)R(y) | |L(x).y+x.R(y)-L(x).R(y) ∪ R(x).y+x.L(x)-R(x).L(x)}.

This definition presupposes the definition of addition.

5. Inverse is defined by

$$x^{-1} = \left\{ 0, \frac{1+\big(R(x)-x\big)L(x)}{R(x)} \cup \frac{1+\big(L(x)-x\big).R(x)}{L(x)} \middle| \frac{1+\big(L(x)-x\big)L(y)}{L(y)} \cup \frac{1+\big(R(x)-x\big)R(y)}{R(y)} \right\}$$

As it is proved in [Conway J.H. 1976] Ch0, 1 the set No is a linearly ordered c-field. The characteristic of No is easily proved to be $\Omega_1$, we call this c-field, c-field of surreal numbers. According to Definition 3 No is an $H_\xi$-leveled Dedekind complete field.

## § 3  The unification .

In this paragraph we prove that all the three different techniques and hierarchies of transfinite real ,of surreal ,of ordinal real numbers  give by      inductive      limit      or      union      the      same      class of numbers .We have already proved that      $CR=\Omega_1 R=C^*R$.      (see corollary 10) and it remains to prove No=CR.

**Lemma 5** . *It holds that $CR=\Omega_1 R=C^*R \subseteq No$.*

*Proof* .Let an open full binary tree $\overset{\bullet}{D}\alpha$ of height      the      principal ordinal a .Then  $\overset{\bullet}{D}\alpha \subseteq No$ ,and the field-inherited operations in the initial segment W($\alpha$) are the Hessenberg operations (see [Conway J.H. 1976] ch 2 § ""containment of the ordinals "note pp 28 and also [Kyritsis C.1991 Alt] the characterisation theorem ).If  $\alpha$ was not a principal ordinal, the W($\alpha$) would not be closed to the Hessenberg operations .Thus the $N_\alpha$, $Z_\alpha$, $Q_\alpha$ are contained in No ,since what it is used to define them from W($\alpha$) is only                                        the field operations .The $Q_\alpha$ is a field  and  from  the    fact    that    No is closed to extensions of its set-subfields (see[ Conway J.H. 1976] ch 4 theorem 28 )we deduce that the field of                     ordinal real numbers $R_\alpha$ is contained in No,  for  every   principal   ordinal number α .Thus $\cup R_\alpha = \Omega_1 R \subseteq No$ .Q.E.D.

**Lemma  6** . *For  every  regular  open  full binary tree  $\overset{\bullet}{D}\alpha$,  it holds that  $\overset{\bullet}{D}\alpha \subseteq R_\beta$,    for  some sufficiently big principal     ordinal number β . (With the inclusion is meant that the restriction of ordering of $R_\alpha$ in the tree, coincides with the ordering of the tree).*

*Proof* . We shall prove it by transfinite induction .It holds  for the trees of finite height. The transfinite induction shall be on the transfinite sequence of all ordinal numbers such that $\aleph(\alpha) = \aleph_{cf(\alpha)}$ and $\aleph(\alpha)$ is a regular aleph. Let us suppose that it holds for all such ordinal numbers of W($\alpha$), and $\aleph(\alpha) = \aleph_{cf(\alpha)}$ and $\aleph(\alpha)$      is      a      regular      aleph      .Then

$: \overset{\bullet}{D}\alpha \subseteq \underset{s \in W(\alpha)}{\cup} \overset{\bullet}{D}_s \underset{s \in W(\alpha)}{\cup} R_{\beta(s)} \subseteq R_{\beta(\alpha)}$ hence $\overset{\bullet}{D}_\alpha \subseteq R_{\beta(\alpha)}$ where $\beta(\alpha)$ is a principal ordinal with $\beta(\alpha) > \underset{s \in W(\alpha)}{\lim} \beta(s)$ . Q.E.D.

From the previous lemma we get that $\cup \overset{\bullet}{D}\alpha = \text{No} \subseteq \Omega_1 R$ ,thus :

**The unification theorem 7**

*It holds that the classes of transfinite real numbers CR , of surreal numbers No, of ordinal real numbers $\Omega_1 R$ ,coincide ,and it is the smallest class (and linearly ordered c-field ) ,that contains all linearly ordered set-fields as subfields.*

We can have obviously analogous statements for the other classes of numbers (complex , quaternion e.t.c.). After the previous theorem, the binary arithnetisation of the order-types, stated in [ Kyritsis C. 1991] II ,theorem 11, is directly provable. We remark that because the levels of the open full binary trees have the property that any upper (lower bounded set has supremum (infimum ) ,(see [Kuratowski K. -Mostowski A 1968] ch ix §1, § 2 theorem 2 ),and after the Hilbert and fundamental (Cauchy) completness of the ordinal real numbers, and remark after definition 13 and ω-normal form according to [ Frankel A.A. 1953] ch 3 theorem 21 ,and after corollary 21 in [Kyritsis C. 1991] ,II , we also get:

**Theorem 8** . *The class of numbers CR=$\Omega_1 R$=No has leveled formal power series representation, leveled Hilbert completeness, leveled fundamental (Cauchy) completeness, leveled $H_\xi$ Dedekind completeness ,leveled supremum completeness and representation with ω-normal forms.*

corps commutatifs de caracteristique zero
C.R. Acad. Sc. Paris, 5 Avril 1961 t. 282
pp 2053-2055.

## List of special symbols

$\alpha, \beta, \omega$          :  Small Greek letters

$\Omega_1$          :  Capital Greek letter omega with  the subscript 1

$F^a$          :  Capital  letter F with  superscript a.

N     :  Capital Aleph ,the first letter of the hebrew alphabet . In the text
is used a capital script. letter n .

$\oplus, \circ$  : cross in a circle, point in a circle .

$N\alpha, Z\alpha, Q\alpha, R\alpha$,:  Roman  capital letters with  subscript  small  Greek  letters

$C\alpha, H\alpha$

$^*X$, $^*R$ et.c  :  Capital standard or roman letters with left superscript a
star.

CN,CZ,CQ,     :Capital  standard  letter  c  followed  by  capital  letters,
with possibly

$C^*R$                a  left superscript  a  star

$\hat{X}$          :  Capital tstandard letter with a cap.

$\Sigma$          :  Capital Greek letter sigma

$\overset{\bullet}{D}_{\alpha}$          :  Capital standard  D with subscript a  small Greek  letter and in
upper place a small zero.